\documentclass[a4,12pt]{article}
\usepackage{amsmath}
\usepackage{graphics}
\usepackage{latexsym}
\usepackage{amsfonts}
\numberwithin{equation}{section} \rightmargin 1.5cm \leftmargin
-1.5cm \textheight 23cm \textwidth 16cm \topmargin -1cm
\title{
 Zero-sum stochastic differential game in finite horizon  involving impulse controls}
\author{Brahim EL ASRI \thanks{Universit\'e Ibn Zohr, Equipe. Aide à la decision,
ENSA, B.P.  1136, Agadir, Maroc. e-mail: b.elasri@uiz.ac.ma }\,\,\,
\, and \, Sehail MAZID \thanks{Universit\'e Ibn Zohr, Equipe. Aide à la decision,
ENSA, B.P.  1136, Agadir, Maroc. e-mail: sehail.mazid@edu.uiz.ac.ma.} }
\begin{document}
\date{}
\maketitle
\newtheorem{theo}{Theorem}
\newtheorem{problem}{Problem}
\newtheorem{pro}{Proposition}
\newtheorem{cor}{Corollary}
\newtheorem{axiom}{Definition}
\newtheorem{rem}{Remark}
\newtheorem{lem}{Lemma}
\newcommand{\brm}{\begin{rem}}
\newcommand{\erm}{\end{rem}}
\newcommand{\beth}{\begin{theo}}
\newcommand{\eeth}{\end{theo}}
\newcommand{\bl}{\begin{lem}}
\newcommand{\el}{\end{lem}}
\newcommand{\bp}{\begin{pro}}
\newcommand{\ep}{\end{pro}}
\newcommand{\bcor}{\begin{cor}}
\newcommand{\ecor}{\end{cor}}
\newcommand{\be}{\begin{equation}}
\newcommand{\ee}{\end{equation}}
\newcommand{\beq}{\begin{eqnarray*}}
\newcommand{\eeq}{\end{eqnarray*}}
\newcommand{\beqa}{\begin{eqnarray}}
\newcommand{\eeqa}{\end{eqnarray}}
\newcommand{\dg}{\displaystyle \delta}
\newcommand{\cm}{\cal M}
\newcommand{\cF}{{\cal F}}
\newcommand{\cR}{{\cal R}}
\newcommand{\bF}{{\bf F}}
\newcommand{\tg}{\displaystyle \theta}
\newcommand{\w}{\displaystyle \omega}
\newcommand{\W}{\displaystyle \Omega}
\newcommand{\vp}{\displaystyle \varphi}
\newcommand{\ig}[2]{\displaystyle \int_{#1}^{#2}}
\newcommand{\integ}[2]{\displaystyle \int_{#1}^{#2}}
\newcommand{\produit}[2]{\displaystyle \prod_{#1}^{#2}}
\newcommand{\somme}[2]{\displaystyle \sum_{#1}^{#2}}
\newlength{\inter}
\setlength{\inter}{\baselineskip} \setlength{\baselineskip}{7mm}
\newcommand{\no}{\noindent}
\newcommand{\rw}{\rightarrow}
\def \ind{1\!\!1}
\def \R{I\!\!R}
\def \N{I\!\!N}
\def \cadlag {{c\`adl\`ag}~}
\def \esssup {\mbox{ess sup}}
\begin{abstract}
This paper considers the problem of two-player zero-sum stochastic differential game with both players adopting impulse controls in finite horizon under  rather  weak  assumptions on the cost functions ($c$ and $\chi$ not decreasing  in time). We use the dynamic programming principle and   viscosity solutions approach to show existence and uniqueness of a solution for  the Hamilton-Jacobi-Bellman-Isaacs (HJBI) partial differential equation (PDE) of the game. We prove that the upper and lower value functions coincide.

\end{abstract}

\no{\bf AMS  subject classifications}:  93E20, 49L20, 49L25, 49N70
\medskip

\no {$\bf Keywords$}:  stochastic differential game, impulse control, quasi-variational inequality, viscosity solution

\section {Introduction}
\no

In this paper  we consider the state process of the stochastic differential game, defined as the solution of the following stochastic equation:

\begin{equation*}
\label{az}
\begin{array}{ll}
 X_{s}=x+\integ{t}{s}b(r,X_{r})dr+\integ{t}{s}\sigma(r,X_{r})dW_{r}
+\sum\limits_{m\geq 1}\xi_{m}\ind_{[\tau_{m},T]}(s)\prod_{l\geq 1}\ind_{\{\tau_{m}\ne\rho_{l}\}}\\ \qquad\qquad\qquad\qquad\qquad\qquad\qquad+\sum\limits_{l\geq 1}\eta_{l}\ind_{[\rho_{l},T]}(s),\qquad\qquad s\geq t,
\end{array}
\end{equation*}
for all $s \in [t,T]$, P-a.s., with $X_{t^-} = x$. Here $W$ is a d-dimensional Wiener process, while
$$u(s)=\sum\limits_{m\geq 1}\xi_{m}\ind_{[\tau_{m},T]}(s)\qquad\qquad\mbox{and}\qquad\qquad v(s)=\sum\limits_{l\geq 1}\eta_{l}\ind_{[\rho_{l},T]}(s)$$
are the impulse controls of player I and player II, respectively. The random variables $\xi_m$ and $\eta_l$ take values in two convex cones $\mathcal{U}$ and $\mathcal{V}$ of $\mathbb{R}^n$, respectively, called the spaces of control actions. The infinite product $\prod_{l\geq 1}\ind_{\{\tau_{m}\ne\rho_{l}\}}$ has the following meaning: When the two players act together on the system at the same time, we take into account only the action of player II. We denote by $X^{t,x,u,v}=\{X^{t,x,u,v}_{s},t\leq s\leq T\}$ the state trajectory of the game with initial time $t$, initial state $x$, and impulse controls $u$ and $v$.
 The gain
functional for player I (resp., cost functional for player II) is
given by
\begin{equation}
\label{aar}
\begin{array}{ll}
J(t, x; u, v)=\mathbb{E}\bigg[\integ{t}{T}f(s,X_{s}^{t,x,u,v})ds-
\sum\limits_{m\geq 1} c(\tau_{m},\xi_{m}) \ind_{[\tau_{m}\leq T]}\prod_{l\geq 1}\ind_{\{\tau_{m}\ne\rho_{l}\}}\\\qquad\qquad\qquad\qquad+\sum\limits_{l\geq 1} \chi(\rho_{l},\eta_{l}) \ind_{[\rho_{l}\leq T]}+g(X_{T}^{t,x,u,v})\bigg],
\end{array}
\end{equation}
$f$ is the running gain and $g$ is the payoff. The function $c$  is the cost function for player I and is a gain function for player II, meaning that when player I performs an action he/she has to pay a cost, resulting in a gain for player II. Analogously, $\chi$ is the cost function for player II and is a gain function for player I. The zero-sum stochastic differential games  problems we will investigate is to  show the upper and lower value functions coincide and the game admits a value.

The stochastic differential games  problem have recently attracted a lot of research activities, especially in connection with mathematical finance, commodities, and in particular energy markets, etc (see e.g. \cite{[TB],[AC],[ACF],[JE],[MJS],[MJP],[RKP],[RKS],[VL],[AJ],[GM],[BO]} and the references therein). In order to tackle those problems, authors use
mainly two approaches. Either a probabilistic one \cite{[JCI],[SHM]} or an approach which uses the Hamilton-Jacobi-Bellman-Isaacs (HJBI) partial differential equation (PDE)\cite{[IKJ]}.

In the finite horizon framework, Cosso \cite{[ACS]} have studied a two-player zero-sum stochastic differential game with both players adopting impulse controls using viscosity solutions to the HJBI equation. They have proved existence and uniqueness for the HJBI equation under stronger constraint (see also Tang and Yong \cite{[ST]}): the cost functions $c$ and $\chi$ are  decreasing in time: $$c(t,y)\geq c(t',y)\qquad \mbox{and}\qquad \chi(t,z)\geq \chi(t',z),$$ for every $0\leq t \leq t'\leq T$, $y\in \mathcal{U}$ and $z\in \mathcal{V}$.

Our aim in this work lies in the fact that we investigate the solution to the zero-sum stochastic differential games  under  rather  weak  assumptions on the cost functions($c$ and $\chi$ are  not decreasing  in time). Therefore the main objective of our work, and this is the novelty of the paper is to characterize the value function as the only solution in viscosity sense of the associated  Hamilton-Jacobi-Bellman-Isaacs (HJBI) partial differential equation (PDE) for the finite horizon problem.

We prove the lower and the upper value functions of stochastic differential game satisfy the dynamic programming principle. We show that the HJBI equation associated to the stochastic differential game, which turns out to be the same for the two value functions because the two players cannot act simultaneously on the system, is the unique solution of the following system:
\begin{equation*} \label{eq:HJBI}\left\{
\begin{array}{c}
max\Big\{ min\Big[-\cfrac{\partial V}{\partial t}-\mathcal{L}V-f,V-H^c_{sup} V\Big],V-H^\chi_{inf} V \Big\}=0
\qquad [0,T)\times\mathbb{R}^n \\ \\
V(T,x)=g(x)\qquad\qquad\qquad\qquad\qquad\qquad\qquad\qquad\qquad\qquad  \forall x\in\mathbb{R}^n.
\end{array}
\right.
\end{equation*}
Where $\mathcal{L}$ is the second-order local operator, and the nonlocal operators $H_{\sup}^c$ and $H_{\inf}^\chi$ are given by
$$H^c_{sup}V(t,x)=\sup\limits_{\xi\in \mathcal{U}}[V(t,x+\xi)-c(t,\xi)],\hspace*{1cm} H^\chi_{inf}V(t,x)=\inf\limits_{\eta\in \mathcal{V}}[V(t,x+\eta)+\chi(t,\eta)],$$
for every $(t,x)\in[0,T]\times\mathbb{R}^n.$

This paper is organized as follows: In Section 2, we formulate the problem and we give the related definitions. In Section 3, we shall introduce the stochastic differential game problem and give some preliminary results of the lower and the upper value functions of stochastic differential game. Further
we provide some estimate for the optimal strategy of problem which in combination with the dynamic programming principle plays a crucial role in the proof of existence the value functions. In Section 4, we prove the dynamic programming principle. Further, we show the existence and continuity of the lower and the upper value functions.  Section 5, is devoted to the connection between the zero sum stochastic differential game problem and Hamilton-Jacobi-Bellman-Isaacs equation. In Section 6, we show that the solution of HJBI is unique in the subclass of bounded continuous functions. Further,  the upper and the lower value functions coincide and the game admits a value.\newpage
\section{Assumptions and formulation of the problem}
\no

Throughout this paper $T$ (resp. $n,\, d$) is a fixed real (resp.
integers) positive numbers. \newline Let us assume the following assumptions:  \medskip

\no(\textbf{H1}) $b:[0,T]\times \mathbb{R}^{n}\rightarrow
\mathbb{R}^{n}$ and $\sigma :[0,T]\times \mathbb{R}^{n}\rightarrow
\mathbb{R}^{n\times d}$ be two continuous functions for which there
exists a constant $C>0$ such that for any $t\in \lbrack 0,T]$ and
$x,x^{\prime }\in \mathbb{R}^{n}$
\begin{equation}
|\sigma (t,x)-\sigma (t,x^{\prime })|+|b(t,x)-b(t,x^{\prime })|\leq
C|x-x^{\prime }|.
 \label{eqs}
\end{equation}%
 Also there exists a constant $C>0$ such that for any $(t,x)\in[0,T]\times\mathbb{R}^n$
\begin{equation}
|\sigma (t,x)|+|b(t,x)|\leq C(1+|x|).
\end{equation}
(\textbf{H2})  $f:[0,T]\times\mathbb{R}^{n} \rightarrow \mathbb{R}$
is uniformly continuous and bounded on $[0,T]\times\mathbb{R}^{n}$. $g:\mathbb{R}^{n} \rightarrow\mathbb{R}$ is uniformly continuous and bounded on $\mathbb{R}^n$.\\
(\textbf{H3}) The cost functions $c:[0,T]\times \mathcal{U} \rightarrow \mathbb{R}$ and $\chi:[0,T]\times \mathcal{V} \rightarrow \mathbb{R}$ are measurable
and uniformly continuous. Furthermore
\begin{equation}\label{cout} \underset{[0,T]\times \mathcal{U}}{inf}c \geq k, \hspace*{2cm} \underset{[0,T]\times \mathcal{V}}{inf}\chi\geq k,
\end{equation}
where $k>0$.\\
Moreover,
\begin{equation} \label{cout1}
  c(t,\xi_1+\xi_2)\leq c(t,\xi_1)+c(t,\xi_2)\\
\end{equation}
\begin{equation}\label{cout2}
  \chi(t,\eta_1+\eta_2)\leq \chi(t,\eta_1)+\chi(t,\eta_2)
\end{equation}
for every $t\in [0,T], \xi_1,\xi_2\in \mathcal{U}$ and $\eta_1,\eta_2\in\mathcal{V}$.
\\
(\textbf{H4}) (no terminal impulse). For any $x\in\mathbb{R}^n, \eta \in\mathcal{V} $ and $\xi\in\mathcal{U}$
\begin{equation}\label{Terminal}
\sup\limits_{\xi\in \mathcal{U}}[g(x+\xi)-c(T,\xi)]\leq g(x) \leq \inf\limits_{\eta\in \mathcal{V}}[g(x+\eta)+\chi(T,\eta)].
\end{equation}
$\mathcal{U}$ and $\mathcal{V}$ are two convex  cones of $\mathbb{R}^{n}$ with $\mathcal{U} \subset \mathcal{V}$.
\begin{rem}
The above assumptions  (\ref{cout1}) and (\ref{cout2})
ensures that multiple impulses occurring at the same time are
suboptimal. $\Box$
\end{rem}

We now consider the HJBI equation:
\begin{equation} \label{eq:HJBI}\left\{
\begin{array}{c}
max\Big\{ min\Big[-\cfrac{\partial V}{\partial t}-\mathcal{L}V-f,V-H^c_{sup} V\Big],V-H^\chi_{inf} V \Big\}=0
\qquad [0,T)\times\mathbb{R}^n, \\ \\
V(T,x)=g(x)\qquad\qquad\qquad\qquad\qquad\qquad\qquad\qquad\qquad\qquad  \forall x\in\mathbb{R}^n.
\end{array}
\right.
\end{equation}
Where $\mathcal{L}$ is the second-order local operator
   $$\mathcal{L}V= \langle b,\nabla _x V\rangle +\cfrac{1}{2}tr[\sigma\sigma^*\nabla^2_xV],$$
and the nonlocal operators $H_{\sup}^c$ and $H_{\inf}^\chi$ are given by
$$H^c_{sup}V(t,x)=\sup\limits_{\xi\in \mathcal{U}}[V(t,x+\xi)-c(t,\xi)],\hspace*{1cm} H^\chi_{inf}V(t,x)=\inf\limits_{\eta\in \mathcal{V}}[V(t,x+\eta)+\chi(t,\eta)],$$
for every $(t,x)\in[0,T]\times\mathbb{R}^n.$

The main objective of this paper is to focus on the existence and uniqueness of
the solution in viscosity sense of (\ref{eq:HJBI}) whose definition
is:
\begin{axiom}
Let $V$ be a continuous function defined on $[0,T]\times\mathbb{R}^n$  and such
that $V(T,x) = g(x)$ for any $x\in \mathbb{R}^n$. The $V$ is called:

(i) A viscosity subsolution of (\ref{eq:HJBI}) if for any $(t_0,x_0)\in[0,T)\times\mathbb{R}^n$ and  any function
 $\phi \in C^{1,2}([0,T)\times\mathbb{R}^n)$, such that $(t_0,x_0)$ is a local maximum  of $V-\phi$, we have:
\begin{equation}
\begin{array}{ll}
max\Big\{ min\Big[-\cfrac{\partial \phi}{\partial t}(t_0,x_0)-\mathcal{L}\phi(t_0,x_0)-f(t_0,x_0),V(t_0,x_0)-H^c_{sup} V(t_0,x_0)\Big],\\\qquad\qquad\qquad V(t_0,x_0)-H^\chi_{inf} V(t_0,x_0) \Big\}\leq 0.
\end{array}
\end{equation}

(ii) A viscosity supersolution of (\ref{eq:HJBI})  if for any $(t_0,x_0)\in[0,T)\times\mathbb{R}^n$ and any function
 $\phi \in C^{1,2}([0,T)\times\mathbb{R}^n)$, such that  is  $(t_0,x_0)$  a local
minimum of $V-\phi$, we have:
\begin{equation}
\begin{array}{ll}
max\Big\{ min\Big[-\cfrac{\partial \phi}{\partial t}(t_0,x_0)-\mathcal{L}\phi(t_0,x_0)-f(t_0,x_0),V(t_0,x_0)-H^c_{sup} V(t_0,x_0)\Big],\\\qquad\qquad\qquad V(t_0,x_0)-H^\chi_{inf} V(t_0,x_0) \Big\}\geq 0.
\end{array}
\end{equation}

(iii) A viscosity solution if it is both a viscosity supersolution and subsolution. $\Box$
\end{axiom}

There is an equivalent formulation of this definition (see e.g. \cite{[MC]}) which we give because it will be
useful later. So firstly, we define the notions of superjet and subjet of a continuous function $V$.
\begin{axiom}
Let $V \in C((0,T)\times\mathbb{R}^n), (t,x)$ an element of $(0,T)\times\mathbb{R}^n$ and finally $S_n$ the set of $n\times n$ symmetric matrices. We denote by $J^{2,+}V(t,x)$ (resp. $J^{2,-}V(t,x))$, the superjets (resp. the subjets) of $\,V$
at $(t,x)$, the set of triples $(p,q,X)\in\mathbb{R}\times\mathbb{R}^n\times S_n$ such that:
\begin{equation*}
\begin{array}{lll}
V(s,y)\leq V(t,x)+p(s-t)+\langle q,y-x\rangle+\displaystyle\frac{1}{2}\langle X(y-x),y-x\rangle+o(|s-t|+|y-x|^2)
\end{array}
\end{equation*}
\begin{equation*}
\begin{array}{lll}
\big(\text{resp.}\\ V(s,y)\geq  V(t,x)+p(s-t)+\langle q,y-x\rangle+\displaystyle\frac{1}{2}\langle X(y-x),y-x\rangle+o(|s-t|+|y-x|^2)\big).\Box
\end{array}
\end{equation*}
\end{axiom}
Note that if $\phi-V$ has a local maximum (resp. minimum) at $(t, x)$, then we obviously have:
$$(D_t\phi(t,x),D_x\phi(t,x),D^2_{xx}\phi(t,x))\in J^{2,-}V(t,x)\;(\text{resp}.\; J^{2,+}V(t,x)).\Box$$

We now give an equivalent definition of a viscosity solution of HJBI equation(\ref{eq:HJBI}):
\begin{axiom}
Let $V$ be a continuous function defined on $[0,T]\times\mathbb{R}^n$  and such
that $V(T,x) = g(x)$ for any $x\in \mathbb{R}^n$. Then V is a viscosity supersolution (resp., subsolution) to the HJBI equation
(\ref{eq:HJBI}) if and only if for every $(t, x)\in [0,T)\times\mathbb{R}^n$ and $(p,q,X)\in J^{2,-}V(t,x)\, ($resp. $J^{2,+}V(t,x))$,
\begin{equation*}
\begin{array}{ll}
max\Big\{ min\Big[-p-\langle b,q\rangle-\displaystyle\frac{1}{2}Tr\big[\sigma^*X\sigma\big] -f(t,x),V(t,x)-H^c_{sup} V(t,x)\Big],\\\qquad\qquad\qquad V(t,x)-H^\chi_{inf} V(t,x)\Big\}\geq 0 \qquad (\text{resp.,}\leq0).
\end{array}
\end{equation*}
\end{axiom}
It is called a viscosity solution it is both a viscosity subsolution and supersolution. $\Box$

As pointed out previously we will show that system (\ref{eq:HJBI})
has a unique solution in viscosity sense. This system is the
deterministic version of the stochastic differential game problem
will describe briefly in the next section.
\section{The stochastic differential game problem}
\subsection{Setting of the problem}
\no

Let ($\Omega,\mathcal{F},\mathbb{P})$ is a fixed probability space on which is defined a
standard $d$-dimensional Brownian motion $W=(W_t)_{t\leq T}$, whose natural filtration is
$(\mathcal{F}^0_{t}:=\sigma\{W_s;s\leq t\})_{0\leq t\leq T}.$ We denote by $\mathbb{F} =(\mathcal{F}_{t})_{t\leq T}$ the completed filtration of $(\mathcal{F}^0_{t})_{t\leq T}$ with the $\mathbb{P}$-null sets of $\mathcal{F}$.
 We are given two convex cones  $\mathcal{U}$ and $\mathcal{V}$ of $\mathbb{R}^n$, with $\mathcal{V}\subset \mathcal{U}$. We call $\mathcal{U}$ and $\mathcal{V}$ the spaces of controle actions. We begin by introducing the concept of impulse control.
\begin{axiom}
An impulse control $u=\sum\limits_{m\geq1}\xi_{m}\ind_{[\tau_{m},T]}$ for player I (resp, $v=\sum\limits_{l\geq1}\eta_{l}\ind_{[\rho_{l},T]}$ for player II) on $[t,T]\subset\mathbb{R}^{+}=[0,+\infty)$, is such that:

(i) $(\tau_{m})_{m}$ (resp., $(\rho_{l})_{l})$, the action times, is a sequence of $\mathbb{F}$-stopping
times, valued in $[t,T]\cup \{+\infty\}$ such that $\mathbb{P}$-a.s. $\tau_{m}\leq \tau_{m+1}$ (resp.,\, $\rho_{l}\leq \rho_{l+1})$.

(ii) $(\xi_{m})_{m}$ (resp., $(\eta_{l})_{l})$, the actions, is a sequence of $\mathcal{U}$-valued (resp, $\mathcal{V}$-valued) random variables, where each $\xi_{m}$ (resp., $\eta_{l}$) is $\mathcal{F}_{\tau_{m}}$-measurable (resp., $\mathcal{F}_{\rho_{l}}$-measurable).
\end{axiom}

Let $t\in[0,T]$ be the initial time of the game and $x \in \mathbb{R}^n$ the initial state. Then, given the impulse controls $u$ and $v$ on $[t,T]$, the state process of the stochastic differential game is defined as the solution to the following stochastic equation:

\begin{equation}
\label{az}
\begin{array}{ll}
 X_{s}=x+\integ{t}{s}b(r,X_{r})dr+\integ{t}{s}\sigma(r,X_{r})dW_{r}
+\sum\limits_{m\geq 1}\xi_{m}\ind_{[\tau_{m},T]}(s)\prod_{l\geq 1}\ind_{\{\tau_{m}\ne\rho_{l}\}}\\ \qquad\qquad\qquad\qquad\qquad\qquad\qquad+\sum\limits_{l\geq 1}\eta_{l}\ind_{[\rho_{l},T]}(s),\qquad\qquad s\geq t,
\end{array}
\end{equation}
Then, under (\textbf{H1}), the stochastic differential equation (\ref{az}) admits a unique solution $X^{t,x,u,v}=\{X^{t,x,u,v}_{s},t\leq s\leq T\}$ for every $(t,x)\in$[0,T]$\times \mathbb{R}^n, u\in\mathcal{U}$ and $v\in\mathcal{V}$. The gain
functional for player I (resp., cost functional for player II) is
given by
\begin{equation}
\label{aar}
\begin{array}{ll}
J(t,x,u,v):=\mathbb{E}\bigg[\integ{t}{T}f(s,X_{s}^{t,x,u,v})ds-
\sum\limits_{m\geq 1} c(\tau_{m},\xi_{m}) \ind_{[\tau_{m}\leq T]}\prod_{l\geq 1}\ind_{\{\tau_{m}\ne\rho_{l}\}}\\\qquad\qquad\qquad\qquad+\sum\limits_{l\geq 1} \chi(\rho_{l},\eta_{l}) \ind_{[\rho_{l}\leq T]}+g(X_{T}^{t,x,u,v})\bigg],
\end{array}
\end{equation}
$f$ is the running gain and $g$ is the payoff. The function $c$  is the cost function for player I and is a gain function for player II, meaning that when player I performs an action he/she has to pay a cost, resulting in a gain for player II. Analogously, $\chi$ is the cost function for player II and is a gain function for player I.

\begin{axiom}
\label{rest}
Let $u=\sum\limits_{m\geq1}\xi_{m} \ind_{[\tau_{m},T]}$  be an  impulse control  on $[t,T]$, and let $\tau \leq \sigma $ be two $[t,T]$-valued $\mathbb{F}$-stopping times. Then we define the restriction $u_{[\tau,\sigma]}$ of the impulse control $u$ by:
\begin{equation} u_{[\tau,\sigma]}(s)=\sum\limits_{m\geq1}\xi_{\mu_{t,\tau}(u)+m}\ind_{[\tau_{\mu_{t,\tau}(u)+m}\leq s\leq \sigma]}(s), \quad\quad\quad \tau \leq s \leq \sigma,
\end{equation}
$\mu_{t,\tau}$ is the number of impulses up to time $\tau$, i.e.,
\begin{equation}\mu_{t,\tau}(u):=\sum\limits_{m\geq 1}\ind_{[\tau_m\leq \tau]}.\end{equation}
\end{axiom}
\begin{axiom}(Admissible impulse control).
An admissible impulse control u for player I (resp., $v$ for player II)
on $[t,T]\subset \mathbb{R}^{+},$ is an impulse control for player I (resp., II) on $[t, T]$ with a finite average number of impulses, i.e.,
$$ \mathbb{E}[\mu_{t,T}(u)]<\infty \hspace*{2cm}(resp.,\,\mathbb{E}[\mu_{t,T}(v)]<\infty).$$
The set of all admissible impulse controls for player I
(resp., II) on $[t,T]$ is denoted by $\mathcal{U}_{t,T}$ (resp., $\mathcal{V}_{t,T})$.
\end{axiom}

Next, we adopt the notion of Control identification and nonanticipative strategy.
\begin{axiom}(Control identification).
\\Let $u=\sum\limits_{m\geq 1}\xi_{m}\ind_{[\tau_{m},T]}$ and $u'=\sum\limits_{m\geq 1}\xi'_{m}\ind_{[\tau'_{m},T]}$
in $\mathcal{U}_{t,T},$ we write $u\equiv u'$ on $[t,T]$
if
$$\mathbb{P}(\{u=u' \;\text{a.e. on}\,\,\, [t,T]\})= 1.$$
Similarly, we interpret $v \equiv v'$ on $[t,T]$ in $\mathcal{V}_{t,T}.$
\end{axiom}
\begin{axiom}(nonanticipative strategy).
The nonanticipative strategy set $\mathcal{A}_{t,T}$ for player I
is the collection of all nonanticipative maps $\alpha$ from  $\mathcal{V}_{t,T}$ to $\mathcal{U}_{t,T}$, i.e  for any $[t,T]$-valued $\mathbb{F}$-stopping times $\tau$ and any $v_{1}$,$v_{2}\in \mathcal{V}_{t,T}$,
\begin{equation*}
\begin{array}{ll}
\text{if}\quad v_{1}\equiv v_{2}\quad \text{on}\quad [t,\tau],\quad \text{then}\\ \\
 \;\alpha(v_{1})\equiv \alpha(v_{2})\quad \text{on}\quad [t,\tau].
\end{array}
\end{equation*}
Analogously, the nonanticipative strategy set $\mathcal{B}_{t,T}$ for player II
is the collection of all nonanticipative maps $\beta$ from  $\mathcal{U}_{t,T}$ to $\mathcal{V}_{t,T}$.
\end{axiom}

We are now ready to introduce the upper and lower values of the game:\\
 For every $(t,x)\in[0,T]\times\mathbb{R}^n$ let us set
\begin{equation}
V^{-}(t,x):=\inf\limits_{\beta\in {\mathcal{B}_{t,T}}}\sup\limits_{u\in\mathcal{U}_{t,T}}J(t,x,u,\beta(u))
\end{equation}
and
\begin{equation}
V^{+}(t,x):=\sup\limits_{\alpha\in \mathcal{A}_{t,T}}\inf\limits_{v\in\mathcal{V}_{t,T}}J(t,x,\alpha(v),v)
\end{equation}
The maps $V^{-}$ and $V^{+}$ are called the lower value and the upper value of the game, respectively. The game is said to admit a value if\, $V^-=V^+$.

The HJBI equation associated to the stochastic differential game, which turns out to be the same for the two value functions because the two players cannot act simultaneously on the system, is given by (\ref{eq:HJBI}).

\begin{rem}
The infinite product $\prod_{l\geq 1}\ind_{\{\tau_{m}\ne\rho_{l}\}}$ in (\ref{az}) means that when the two players act together on the system, we take into account only the action of player II. When  take into account only the action of player I instead of player II. Then, using arguments analogous to those presented in the sections below, it can be proved that, with this assumption, the corresponding HJBI equation is given by
\end{rem}
\begin{equation} \label{eq:HJBI1}\left\{
\begin{array}{c}
min\Big\{ max\Big[-\cfrac{\partial V}{\partial t}-\mathcal{L}V-f,V-H^\chi_{inf} V\Big],V-H^c_{sup} V \Big\}=0
\qquad [0,T)\times\mathbb{R}^n \\ \\
V(T,x)=g(x)\qquad\qquad\qquad\qquad\qquad\qquad\qquad\qquad\qquad\qquad  \forall x\in\mathbb{R}^n.
\end{array}
\right.
\end{equation}
\subsection{Preliminary results}
\no

 In this section we present some properties of the lower and upper value functions
of our differential game.

 We begin by the following lemma, which is concerned with the continuous
dependence of $X^{t,x,u,v}$ with respect to $x$.
\begin{lem}
\label{lem:2}
 Under assumption (\textbf{H1}) there exists a constant $C>0$ such that, for every $t\in[0,T],\,x,x'\in\mathbb{R}^n, u\in \mathcal{U}_{t,T}$ and $ v\in \mathcal{V}_{t,T}$ we have:
\begin{equation}
 \mathbb{E}\big[|X_s^{t,x,u,v}-X_s^{t,x',u,v}|\big]\leq C|x-x'|.
\end{equation}
\end{lem}
$Proof$. see [\cite{[HP]}, Appendix]

Next, in the following proposition, we prove that the two value functions are
bounded.
\begin{theo}
Under the standing assumptions (\textbf{H1}), (\textbf{H2}) and (\textbf{H3}), the lower and upper value functions are bounded.
\end{theo}
$Proof$: We make the proof only for the lower value function $V^-$, the other case being
analogous.\\
Using the definition of lower value function, we have, for $(t,x)\in[0,T]\times\mathbb{R}^n$
\begin{multline*}
V^-(t,x)=\inf\limits_{\beta\in {\mathcal{B}_{t,T}}}\sup\limits_{u\in\mathcal{U}_{t,T}}\mathbb{E}\bigg[\integ{t}{T}f(s,X_{s}^{t,x,u,\beta(u)})ds-
\sum\limits_{m\geq 1} c(\tau_{m},\xi_{m}) \ind_{[\tau_{m}\leq T]}\prod_{l\geq 1}\ind_{\{\tau_{m}\ne\rho_{l}\}}\\+\sum\limits_{l\geq 1} \chi(\rho_{l},\eta_{l}) \ind_{[\rho_{l}\leq T]}+g(X_{T}^{t,x,u,\beta(u)})\bigg] \\ \leq \sup\limits_{u\in\mathcal{U}_{t,T}}\mathbb{E}\bigg[\integ{t}{T}f(s,X_{s}^{t,x,u,v_0})ds-\sum\limits_{m\geq 1} c(\tau_{m},\xi_{m}) \ind_{[\tau_{m}\leq T]}+g(X_{T}^{t,x,u,v_0})\bigg],
\end{multline*}
where $v_0$ is the control
with no impulses.\\
Let $\epsilon > 0$, then there exists a strategy $u^\epsilon\in \mathcal{U}_{t,T}$ such that
$$ V^-(t,x)\leq \mathbb{E}\bigg[\integ{t}{T}f(s,X_{s}^{t,x,u^\epsilon,v_0})ds-\sum\limits_{m\geq 1} c(\tau^\epsilon_{m},\xi^\epsilon_{m}) \ind_{[\tau^\epsilon_{m}\leq T]}+g(X_{T}^{t,x,u^\epsilon,v_0})\bigg]+\epsilon.
$$
Since the cost $c(\tau^\epsilon_{m},\xi^\epsilon_{m})$ are non negative functions, then we have
$$ V^-(t,x)\leq \mathbb{E}\bigg[\integ{t}{T}f(s,X_{s}^{t,x,u^\epsilon,v_0})ds+g(X_{T}^{t,x,u^\epsilon,v_0})\bigg]+\epsilon.
$$
Therefore we find, using also the boundedness of $f$ and $g$, that there exists a constant $C>0$ such that
$$ V^{-}(t,x)\leq C.$$
In a similar way we can prove that there exists a constant $C>0$ such that
$$ V^{-}(t,x)\geq -C,$$
from which we deduce the thesis.$\quad\Box$

We are now giving some properties of nearly nearly optimal strategies.
\begin{pro}
\label{pro1}
Let $u\in\mathcal{U}_{t,T}$ and $v\in\mathcal{V}_{t,T}$ be a  nearly optimal strategies composed of impulses control $(\tau,\xi)=((\tau_n)_{n\geq1},(\xi_n)_{n\geq1})$ and $(\rho,\eta)=((\rho_n)_{n\geq1},(\eta_n)_{n\geq1}).$ Then:
\begin{equation}
\label{condition1}
\sum\limits_{m\geq 1} c(\tau_{m},\xi_{m}) \ind_{[\tau_{m}\leq T]}+\sum\limits_{l\geq 1} \chi(\rho_{l},\eta_{l}) \ind_{[\rho_{l}\leq T]}\leq C.
\end{equation}
There exists a positive constant $C>0$ which does not depend on $t$ and $x$ such that:
\begin{equation}
\label{condition2}
\forall n\geq 1,\qquad \mathbb{P}[\tau_n\leq T]+\mathbb{P}[\rho_n\leq T]\leq\frac{C}{n}.
\end{equation}
\end{pro}
We denote by $\widehat{\mathcal{U}}_{t,T}$ and $\widehat{\mathcal{V}}_{t,T}$ the set which satisfies the conditions (\ref{condition1}) and (\ref{condition2}), respectively. Similarly, by $\widehat{\mathcal{A}}_{t,T}$ and $\widehat{\mathcal{B}}_{t,T}$ the sets that contain all the
nonanticipative strategies with values in $\widehat{\mathcal{U}}_{t,T}$ and $\widehat{\mathcal{V}}_{t,T}$, respectively.\\
$Proof$.
Let us choose a nearly optimal strategy $u\in\mathcal{U}_{t,T}$ composed of impulse control $((\tau_n)_{n\geq 1},(\xi_n)_{n\geq 1})$ such that,
for $(t,x)\in[0,T]\times\mathbb{R}^n$,
\begin{multline*} \mathbb{E}\bigg[\integ{t}{T}f(s,X_{s}^{t,x,u,v_0})ds-\sum\limits_{m\geq 1} c(\tau_{m},\xi_{m}) \ind_{[\tau_{m}\leq T]}+g(X_{T}^{t,x,u,v_0})\bigg]
\geq V^-(t,x)-1.
\end{multline*}
Since $V^-$, $f$ and $g$ are bounded, then we have
$$\mathbb{E}\bigg[\sum\limits_{m\geq 1} c(\tau_{m},\xi_{m}) \ind_{[\tau_{m}\leq T]}\bigg]\leq C.$$
Next we show $(\ref{condition2})$. Taking into account that $c(\tau,\xi)\geq k > 0$ for any $(\tau,\xi)\in[0,T]\times\mathcal{U}_{t,T}$, we obtain:
$$\mathbb{E}\bigg[\sum\limits_{m\geq 1} k \ind_{[\tau_{m}\leq T]}\bigg]\leq C.$$ But for any $ m\leq n,\, [\tau_n \leq T ]\subseteq [\tau_m\leq T]$, then:
$$ \mathbb{E}\big[n k\ind_{[\tau_{n}\leq T]}\big]\leq C.$$
Finally taking into account $k > 0$, we obtain the desired result.\\ The other case being analogous.
\quad\quad\qquad\qquad\quad\qquad\qquad\qquad\qquad\qquad\qquad\qquad$\Box$
\begin{cor}\label{cor1}
Under the standing assumptions (\textbf{H1}), (\textbf{H2}) and (\textbf{H3}), $J(t,x,u,v)$ is bounded for every $(t,x)\in [0,T]\times\mathbb{R}^n, u\in\widehat{\mathcal{U}}_{t,T}$ and $ v\in\widehat{\mathcal{V}}_{t,T}$.
\end{cor}

In the following proposition, using Lemma \ref{lem:2}, we prove that the lower and upper
value functions are continuous in the state variable, together with the gain
functional.
\begin{pro}
\label{pro2}
The  gain functional, lower and upper
value functions are continuous in $x$.
\end{pro}
$Proof.$ It is enough to show that the conclusion holds for the gain functional $J$.\\
For every $t\in[0,T],\, x,x'\in \mathbb{R}^n, u\in\mathcal{U}_{t,T}, v\in\mathcal{V}_{t,T}$, we have
\begin{multline*}
\qquad\qquad |J(t,x,u,v)-J(t,x',u,v)|\leq \mathbb{E}\bigg[\int_{t}^{T}|f(s,X_s^{t,x,u,v})-f(s,X_s^{t,x',u,v})|ds\\+|g(X_T^{t,x,u,v})-g(X_T^{t,x',u,v})|\bigg].\qquad\qquad\qquad
\qquad\qquad
\end{multline*}
From Lemma $\ref{lem:2}$ and continuity of $f$ and $g$ in $x$ we get the thesis.\qquad\qquad$\Box$
\section{The value functions}
\subsection{Dynamic Programming Principle}
\no

 In this section we present the dynamic programming principle (DPP) for the stochastic
differential game. We begin with the following lemma.
\begin{lem}(\cite{[ACS]}, Lemma 3.3).
The lower and upper
value functions are given by
\begin{equation}
V^{-}(t,x):=\inf\limits_{\beta\in \bar{{\mathcal{B}}}_{t,T}}\sup\limits_{u\in\bar{\mathcal{U}}_{t,T}}J(t,x,u,\beta(u))
\end{equation}
and
\begin{equation}
V^{+}(t,x):=\sup\limits_{\alpha\in \bar{\mathcal{A}}_{t,T}}\inf\limits_{v\in\bar{\mathcal{V}}_{t,T}}J(t,x,\alpha(v),v)
\end{equation}
for every $(t,x)\in[0,T)\times\mathbb{R}^n$, where $\bar{\mathcal{U}}_{t,T}$ and $\bar{\mathcal{V}}_{t,T}$  contain all the impulse controls
in $\widehat{\mathcal{U}}_{t,T}$ and $\widehat{\mathcal{V}}_{t,T}$, respectively, which have no impulses at time t. Similarly, $\bar{\mathcal{A}}_{t,T}$ and $\bar{{\mathcal{B}}}_{t,T}$ are subsets of $\widehat{\mathcal{A}}_{t,T}$ and $\widehat{\mathcal{B}}_{t,T}$, respectively. In particular, they contain all the nonanticipative strategies with values in $\bar{\mathcal{U}}_{t,T}$ and $\bar{\mathcal{V}}_{t,T}$, respectively.
\end{lem}
\begin{theo}
Under assumptions (\textbf{H1}), (\textbf{H2}) and (\textbf{H3}), given $0\leq t\leq s\leq T$, $x\in \mathbb{R}^n$, and each family of $\mathbb{Q}_{[t,T]}-$valued $\mathcal{F}_{t,s}$-stopping times $\{\tau^{u,v},(u,v)\in\widehat{\mathcal{U}}_{t,T}\times\widehat{\mathcal{V}}_{t,T}\}$, we have
\begin{equation}
\begin{array}{ll}
V^{-}(t,x)\leq\inf\limits_{\beta\in \widehat{\mathcal{B}}_{t,T}}\sup\limits_{u\in\widehat{\mathcal{U}}_{t,T}}\mathbb{E}\bigg[\integ{t}{\tau^{u,\beta}}f(r,X_{r}^{t,x,u,\beta(u)})dr-
\sum\limits_{m\geq 1} c(\tau_{m},\xi_{m}) \ind_{[\tau_{m}\leq \tau^{u,\beta}]}\prod_{l\geq 1}\ind_{\{\tau_{m}\ne\rho_{l}\}}\\\qquad\qquad\qquad\qquad +\sum\limits_{l\geq 1} \chi(\rho_{l},\eta_{l}) \ind_{[\rho_{l} \leq \tau^{u,\beta}]}+(V^{-})^*(\tau^{u,\beta},X_{\tau^{u,\beta}}^{t,x,u,\beta(u)})\bigg],
\end{array}
\end{equation}
\begin{equation}
\begin{array}{ll}
V^{-}(t,x)\geq\inf\limits_{\beta\in \widehat{\mathcal{B}}_{t,T}}\sup\limits_{u\in\widehat{\mathcal{U}}_{t,T}}\mathbb{E}\bigg[\integ{t}{\tau^{u,\beta}}f(r,X_{r}^{t,x,u,\beta(u)})dr-
\sum\limits_{m\geq 1} c(\tau_{m},\xi_{m}) \ind_{[\tau_{m}\leq \tau^{u,\beta}]}\prod_{l\geq 1}\ind_{\{\tau_{m}\ne\rho_{l}\}}\\\qquad\qquad\qquad\qquad +\sum\limits_{l\geq 1} \chi(\rho_{l},\eta_{l}) \ind_{[\rho_{l} \leq \tau^{u,\beta}]}+(V^{-})_*(\tau^{u,\beta},X_{\tau^{u,\beta}}^{t,x,u,\beta(u)})\bigg],
\end{array}
\end{equation}
\begin{equation}
\begin{array}{ll}
V^{+}(t,x)\leq\sup\limits_{\alpha\in\widehat{\mathcal{A}}_{t,T}}\inf\limits_{v\in {\widehat{\mathcal{V}}_{t,T}}}\mathbb{E}\bigg[\integ{t}{\tau^{\alpha,v}}f(r,X_{r}^{t,x,\alpha(v),v})dr-
\sum\limits_{m\geq 1} c(\tau_{m},\xi_{m}) \ind_{[\tau_{m}\leq \tau^{\alpha,v}]}\prod_{l\geq 1}\ind_{\{\tau_{m}\ne\rho_{l}\}}\\\qquad\qquad\qquad\qquad +\sum\limits_{l\geq 1} \chi(\rho_{l},\eta_{l}) \ind_{[\rho_{l}\leq \tau^{\alpha,v}]}+(V^{+})^*(\tau^{\alpha,v},X_{\tau^{\alpha,v}}^{t,x,\alpha(v),v})\bigg],
\end{array}
\end{equation}
and
\begin{equation}
\begin{array}{ll}
V^{+}(t,x)\geq\sup\limits_{\alpha\in\widehat{\mathcal{A}}_{t,T}}\inf\limits_{v\in {\widehat{\mathcal{V}}_{t,T}}}\mathbb{E}\bigg[\integ{t}{\tau^{\alpha,v}}f(r,X_{r}^{t,x,\alpha(v),v})dr-
\sum\limits_{m\geq 1} c(\tau_{m},\xi_{m}) \ind_{[\tau_{m}\leq \tau^{\alpha,v}]}\prod_{l\geq 1}\ind_{\{\tau_{m}\ne\rho_{l}\}}\\\qquad\qquad\qquad\qquad +\sum\limits_{l\geq 1} \chi(\rho_{l},\eta_{l}) \ind_{[\rho_{l}\leq \tau^{\alpha,v}]}+(V^{+})_*(\tau^{\alpha,v},X_{\tau^{\alpha,v}}^{t,x,\alpha(v),v})\bigg],
\end{array}
\end{equation}
where
 $\mathbb{Q}_{[t,T]}=(\mathbb{Q} \cap[t,T])\cup\{t,T\},$ and  $V_*$ (resp., $V^*$) its lower (resp.
upper) semicontinuous envelope defined by:
\begin{equation*}
V_*(t,x):=\liminf\limits_{\substack{(s,y)\to(t,x)\\
(s,y)\in[0,T]\times\mathbb{R}^n }}V(s,y) \qquad\text{and}\qquad V^*(t,x):=\limsup\limits_{\substack{(s,y)\to(t,x)\\
(s,y)\in[0,T]\times\mathbb{R}^n }}V(s,y).
\end{equation*}

\end{theo}
$Proof$. We prove the dynamic programming principle only for $V^-$, the other case
being analogous.

Let $\epsilon > 0$ and consider an arbitrary function
$$\phi:[0,T]\times\mathbb{R}^n \longrightarrow \mathbb{R} \quad \text{such that}\quad \phi\; \text{continuous, bounded from above and}\; V^-\leq\phi.$$
For each $(s,y)\in[0,T]\times\mathbb{R}^n$, there exists $\beta^{\epsilon}\in \bar{\mathcal{B}}_{s,T}$ such that
\begin{equation}
\label{progeq1}
\phi(s,y)\geq V^-(s,y)\geq J(s,y,u_{[s,T]},\beta^{\epsilon}(u_{[s,T]}))-\epsilon,
\end{equation}
where $u_{[s,T]}$ is as introduced in Definition \ref{rest}.
Let $\mathbb{Q}_{[t,T]}:=\{t_i\}_{i\geq1}$ and  fix one of the points $t_i$ in time. For each $y\in\mathbb{R}^n$, the continuity of $J$ established in Proposition $\ref{pro2}$ and that of $\phi$ imply that there exists a neighborhood $B^i(y)\subseteq\mathbb{R}^n$ (of size depending on $y, i, \epsilon$) such that
\be \left\{
\begin{array}{ll}\label{progeq2}
\phi(t_i,y')\geq\phi(t_i,y)-\epsilon,\\
J(t_i,y',u_{[t_i,T]},v)\leq J(t_i,y,u_{[t_i,T]},v)+\epsilon,
\end{array}
\right. \ee
for all $y'\in B^i(y)$ and $v\in\bar{\mathcal{V}}_{t_i,T}$.\\
Therefore, since the family $\{B^i(y):y\in\mathbb{R}^n\}$ forms an open cover of $\mathbb{R}^n$, there exists a sequence $(y_j)_{j\geq1}$ in $\mathbb{R}^n$ such that $\{B^i(y_j)\}_{j\geq1}$ is a countable subcover of $\mathbb{R}^n$. We set $\beta_{i,j}:=\beta^{t_i,y_j}\in\bar{\mathcal{B}}_{t_i,T}$ and $B^i_j:=B^i(y_j)$.\\
We can now define, for $i$ still being fixed, a measurable partition $(A^i_j)_{j\geq1}$  by
$$ A^i_1:=B^i_1,\quad A^i_{j+1}:=B^i_{j+1}\setminus(B^i_1\cup...\cup B^i_j),\quad j\geq1.$$
Since $A^i_j\subseteq B^i_j$, the inequalities $(\ref{progeq1})$ and $(\ref{progeq2})$ yield that
\begin{equation}\label{progeq3}
\phi(t_i,y')\geq J(t_i,y',u_{[t_i,T]},\beta_{i,j}(u_{[t_i,T]})) -3\epsilon \quad\text{for all}\; y'\in A^i_j.
\end{equation}
Now, let $(u,\beta)\in\widehat{\mathcal{U}}_{t,T}\times\widehat{\mathcal{B}}_{t,T}$ be arbitrary and set $\tau=\tau^{\beta(u),u}$. Fix an integer $k\geq1$, we now focus on $(t_i)_{1\leq i\leq k}.$ We
may assume that $t_1<t_2<...<t_k,$ by eliminating and relabeling some of the $t_i$. We define the $\mathcal{F}^t_\tau$-measurable sets
$$ \Gamma^i_j:=\{\tau=t_i\;\text{and}\; X_{t_i}^{t,x,u,\beta(u)}\in A^i_j\}\in \mathcal{F}_{t_i} \quad \text{and} \quad \Gamma(k):=\bigcup_{1\leq i,j\leq k}\Gamma^i_j.$$
Since the $t_i$ are distinct and $A^i_j\cap A^i_{j'}=\emptyset$ for $j\ne j'$, we have $\Gamma^i_j\cap \Gamma^{i'}_{j'}=\emptyset$ for $(i,j)\ne(i',j').$ We construct the strategy $\beta^k$ by
$$ \beta^k(u)=\beta(u)\ind_{[t,\tau]}+\ind_{(\tau,T]}\bigg(\beta(u)\ind_{(\Gamma(k))^c}+\sum\limits_{1\leq i,j\leq k}\beta_{i,j}(u_{[t_j,T]})\ind_{\Gamma^i_j}\bigg),$$
Then
\begin{equation}
\begin{array}{ll}
J(t,x,u,\beta^{k}(u))=  \mathbb{E}\bigg[\integ{t}{\tau}f(r,X_{r}^{t,x,u,\beta(u)})dr-
\sum\limits_{m\geq 1}c(\tau_{m},\xi_{m})\ind_{[\tau_{m}\leq \tau]}\prod_{l\geq 1}\ind_{\{\tau_{m}\ne\rho_{l}\}}\\ \\ +\sum\limits_{l\geq 1} \chi(\rho_{l},\eta_{l})\ind_{[\rho_{l}\leq \tau]}+J(t,x,u,\beta(u))\ind_{(\Gamma(k))^c}+\sum\limits_{1\leq i,j\leq k}J(\tau,X_{\tau}^{t,x,u,\beta(u)},u_{[t_i,T]},\beta_{i,j}(u_{[t_i,T]}))\ind_{\Gamma_j^i}\bigg].
\end{array}
\end{equation}
We deduce via $(\ref{progeq3})$ that
\begin{equation}\label{progeq4}
\mathbb{E}\bigg[\sum\limits_{1\leq i,j\leq k}J(\tau,X_{\tau}^{t,x,u,\beta(u)},u_{[t_i,T]},\beta_{i,j}(u_{[t_i,T]}))\ind_{\Gamma_j^i}\bigg]\leq\mathbb{E}\big[\phi(\tau,X_\tau^{t,x,u,\beta(u)})\ind_{\Gamma(k)}\big]+3\epsilon,
\end{equation}
for every $k\geq1$. Letting $k\rightarrow\infty$, therefore,
$$\mathbb{E}\big[J(t,x,u,\beta(u))\ind_{(\Gamma(k))^c}\big]\longrightarrow 0$$
by dominated convergence and Corollary \ref{cor1}. Moreover, monotone convergence yields
$$\mathbb{E}\big[\phi(\tau,X_\tau^{t,x,u,\beta(u)})\ind_{\Gamma(k)}\big]\longrightarrow \mathbb{E}\big[\phi(\tau,X_\tau^{t,x,u,\beta(u)})\big].$$
Therefore, we deduce the existence of an integer $k_0 \geq 1$ such that
\begin{equation}\label{progeq5}
\begin{array}{ll}
J(t,x,u,\beta^{k_0}(u))\leq  \mathbb{E}\bigg[\integ{t}{\tau}f(r,X_{r}^{t,x,u,\beta(u)})dr-
\sum\limits_{m\geq 1}c(\tau_{m},\xi_{m})\ind_{[\tau_{m}\leq \tau]}\prod_{l\geq 1}\ind_{\{\tau_{m}\ne\rho_{l}\}}\\ +\sum\limits_{l\geq 1} \chi(\rho_{l},\eta_{l})\ind_{[\rho_{l}\leq \tau]}+\phi(\tau,X_\tau^{t,x,u,\beta(u)})\bigg]+4\epsilon.
\end{array}
\end{equation}
Let a sequence of continuous functions $(\phi_n)_n$ such that $\phi_n\geq(V^-)^*$ for all $n\geq1$ and such that $\phi_n$ converges pointwise to $(V^-)^*$.
Set $\varphi_N:=max_{n\geq N}\phi_n$ for $N\geq 1$ and observe that the sequence $(\varphi_N)_N$ is decreasing
and converges pointwise to $(V^-)^*.$ By the
monotone convergence theorem, we then obtain:
$$\mathbb{E}\big[\varphi_N(\tau,X_\tau^{t,x,u,\beta(u)})\big]
\underset{N\to+\infty}{\longrightarrow} \mathbb{E}\big[(V^-)^*(\tau,X_\tau^{t,x,u,\beta(u)})\big].$$
Therefore, we have
\begin{equation}
\begin{array}{ll}
\label{progeq6}
V^-(t,x) \leq  \sup\limits_{u\in\widehat{\mathcal{U}}_{t,T}}\mathbb{E}\bigg[\integ{t}{\tau}f(r,X_{r}^{t,x,u,\beta(u)})dr-
\sum\limits_{m\geq 1} c(\tau_{m},\xi_{m}) \ind_{[\tau_{m}\leq \tau]}\prod_{l\geq 1}\ind_{\{\tau_{m}\ne\rho_{l}\}}\\ \qquad\qquad +\sum\limits_{l\geq 1} \chi(\rho_{l},\eta_{l}) \ind_{[\rho_{l} \leq \tau]}+(V^-)^*(\tau,X_\tau^{t,x,u,\beta(u)})\bigg]+4\epsilon.
\end{array}
\end{equation}
As $\beta$  and $\epsilon$ are arbitrary then sending $\epsilon \to 0,$ and take infimums in right-hand side of $(\ref{progeq6})$, to obtain:
\begin{equation}
\begin{array}{ll}
V^-(t,x)\leq \inf\limits_{\beta\in {\widehat{\mathcal{B}}_{t,T}}}\sup\limits_{u\in\widehat{\mathcal{U}}_{t,T}}\mathbb{E}\bigg[\integ{t}{\tau}f(r,X_{r}^{t,x,u,\beta(u)})dr-
\sum\limits_{m\geq 1} c(\tau_{m},\xi_{m}) \ind_{[\tau_{m}\leq \tau]}\prod_{l\geq 1}\ind_{\{\tau_{m}\ne\rho_{l}\}}\\ \qquad\qquad +\sum\limits_{l\geq 1} \chi(\rho_{l},\eta_{l}) \ind_{[\rho_{l} \leq \tau]}+(V^-)^*(\tau,X_\tau^{t,x,u,\beta(u)})\bigg].
\end{array}
\end{equation}
In a similar way we can prove the reverse inequality, hence deducing the thesis. \qquad\qquad$\Box$
\subsection{Continuity of value functions in time}
\no

 In this section we prove the continuity of\, lower value function and upper value function in t .

 First, we know that any stopping time $\tau$  can be approximated
by a sequence of rational stopping times $\tau^{(n)}=2^{-n}([2^n\tau]-1)$.  Then we can pose

(\textbf{H5}) The action times $\rho_l$ and $\tau_m$ takes values in $\mathbb{Q}_{[t,T]}$
 \beth The lower and upper
value functions are continuous in
$t$.\eeth
$Proof.$ We make the proof only for $V^{-}$, the other case being analogous.\\
First  let us show that $V^-$ is lower semi-continuous.\\
Recall the characterization of dynamical programming principle that reads as
\begin{equation*}
\begin{array}{ll}
V^{-}(t',x)\leq\inf\limits_{\beta\in {\widehat{\mathcal{B}}_{t',T}}}\sup\limits_{u\in\widehat{\mathcal{U}}_{t',T}}\mathbb{E}\bigg[\integ{t'}{\rho_n\wedge T}f(r,X_{r}^{t',x,u,\beta(u)})dr-
\sum\limits_{ m\geq 1} c(\tau_{m},\xi_{m}) \ind_{[\tau_{m}\leq \rho_n\wedge T]}\prod_{l\geq 1}\ind_{\{\tau_{m}\ne\rho_{l}\}}\\ \qquad +\sum\limits_{1\leq l\leq n} \chi(\rho_{l},\eta_{l}) \ind_{[\rho_{l}\leq T]}+\ind_{[\rho_n\leq T]}(V^{-})^*(\rho_n,X_{\rho_n}^{t',x,u,\beta(u)})+\ind_{[\rho_n=+\infty]}g(X_{T}^{t',x,u,\beta(u)})\bigg].
\end{array}
\end{equation*}
\begin{equation*}
\begin{array}{ll}
V^{-}(t,x)\geq\inf\limits_{\beta\in {\widehat{\mathcal{B}}_{t,T}}}\sup\limits_{u\in\widehat{\mathcal{U}}_{t,T}}\mathbb{E}\bigg[\integ{t}{\rho_n\wedge T}f(r,X_{r}^{t,x,u,\beta(u)})dr-
\sum\limits_{ m\geq 1} c(\tau_{m},\xi_{m}) \ind_{[\tau_{m}\leq \rho_n\wedge T]}\prod_{l\geq 1}\ind_{\{\tau_{m}\ne\rho_{l}\}}\\ \qquad +\sum\limits_{1\leq l\leq n} \chi(\rho_{l},\eta_{l}) \ind_{[\rho_{l}\leq T]}+\ind_{[\rho_n\leq T]}(V^{-})_*(\rho_n,X_{\rho_n}^{t,x,u,\beta(u)})+\ind_{[\rho_n=+\infty]}g(X_{T}^{t,x,u,\beta(u)})\bigg].
\end{array}
\end{equation*}
Fix an arbitrary $\epsilon>0$, and  we assume that $t < t'$. Pick $\beta^\epsilon\in\widehat{\mathcal{B}}_{t,T}$ such that
\begin{equation*}
\begin{array}{ll}
V^{-}(t,x)+\epsilon\geq \mathbb{E}\bigg[\integ{0}{\rho^\epsilon_n\wedge T}f(r,X_{r}^{t,x,u,\beta^\epsilon(u)})\ind_{[r\geq t]}dr-
\sum\limits_{ m\geq 1} c(t\vee\tau_{m},\xi_{m}) \ind_{[\tau_{m}\leq \rho_n\wedge T]}\prod_{l\geq 1}\ind_{\{\tau_{m}\ne\rho^\epsilon_{l}\}}\\ \qquad +\sum\limits_{1\leq l\leq n} \chi(t\vee\rho^\epsilon_{l},\eta^\epsilon_{l}) \ind_{[\rho^\epsilon_{l}\leq T]}+\ind_{[\rho^\epsilon_n\leq T]}(V^{-})_*(\rho^\epsilon_n,X_{\rho^\epsilon_n}^{t,x,u,\beta^\epsilon(u)})+\ind_{[\rho^\epsilon_n=+\infty]}g(X_{T}^{t,x,u,\beta^\epsilon(u)})\bigg],\qquad\qquad\qquad\qquad
\end{array}
\end{equation*}
where $u\in\widehat{\mathcal{U}}_{t,T}$ will be chosen later. On the other hand, pick $u^\epsilon\in\widehat{\mathcal{U}}_{t',T}$ such that
\begin{equation*}
\begin{array}{ll}
V^{-}(t',x)-\epsilon\\ \\ \leq\mathbb{E}\bigg[\integ{0}{\rho^\epsilon_n\wedge T}f(r,X_{r}^{t',x,u^\epsilon,\bar{\beta}(u^\epsilon)})\ind_{[r\geq t']}dr-
\sum\limits_{ m\geq 1} c(t'\vee\tau^\epsilon_{m},\xi^\epsilon_{m}) \ind_{[\tau^\epsilon_{m}\leq \rho^\epsilon_n\wedge T]}\prod_{l\geq 1}\ind_{\{\tau^\epsilon_{m}\ne\rho^\epsilon_{l}\}}\\ +\sum\limits_{1\leq l\leq n} \chi(t'\vee\rho^\epsilon_{l},\eta^\epsilon_{l}) \ind_{[\rho^\epsilon_{l}\leq T]}+\ind_{[\rho^\epsilon_n\leq T]}(V^{-})_*(\rho^\epsilon_n,X_{\rho^\epsilon_n}^{t',x,u^\epsilon,\bar{\beta}(u^\epsilon)})+\ind_{[\rho^\epsilon_n=+\infty]}g(X_{T}^{t',x,u^\epsilon,\bar{\beta}(u^\epsilon)})\bigg],\qquad\qquad\qquad\qquad
\end{array}
\end{equation*}
where $$\bar{\beta}(u^\epsilon)=\sum\limits_{1\leq l\leq n} \eta^\epsilon_{l}\ind_{[t'\vee\rho^\epsilon_{l},T]}\in\mathcal{B}_{t',T}.$$
Therefore
\begin{equation*}
\begin{array}{ll}
\label{equ22}
V^{-}(t,x)-V^{-}(t',x) \geq \mathbb{E}\bigg[\integ{0}{\rho_n^\epsilon\wedge T}f(r,X_{r}^{t,x,u^\epsilon,\beta^{\epsilon}(u^\epsilon)})\ind_{[r\geq t]}dr+
\sum\limits_{1\leq l\leq n} \chi(t\vee\rho^\epsilon_{l},\eta^\epsilon_{l})\ind_{[\rho^\epsilon_{l}\leq T]} \\ \\ +\ind_{[\rho_n^\epsilon\leq T]}(V^{-})_*(\rho_n^\epsilon,X_{\rho_n^\epsilon}^{t,x,u^\epsilon,\beta^{\epsilon}(u^\epsilon)})+\ind_{[\rho_n^\epsilon=+\infty]}g(X_{T}^{t,x,u^\epsilon,\beta^\epsilon(u^\epsilon)})\bigg]\\ \\-\mathbb{E}\bigg[\integ{0}{\rho^\epsilon_n\wedge T}f(r,X_{r}^{t',x,u^\epsilon,\bar{\beta}(u^\epsilon)})\ind_{[r\geq t']}dr+
\sum\limits_{1\leq l\leq n}\chi(t'\vee\rho_l^\epsilon,\eta^\epsilon_{l})\ind_{[\rho^\epsilon_l\leq T]}\\ \\+\ind_{[\rho_n^\epsilon\leq T]}(V^{-})^*(\rho_n^\epsilon,X_{\rho_n^\epsilon}^{t',x,u^\epsilon,\bar{\beta}(u^\epsilon)})+\ind_{[\rho_n^\epsilon=+\infty]}g(X_{T}^{t',x,u^\epsilon,\bar{\beta}(u^\epsilon)})\bigg]-2\epsilon\qquad\qquad\qquad
\\ \\ \geq \mathbb{E}\bigg[
-\integ{0}{\rho_n^\epsilon\wedge T}\displaystyle\{|f(r,X_{r}^{t,x,u^{\epsilon},\beta^{\epsilon}({u^{\epsilon}})})-f(r,X_{r}^{t',x,u^\epsilon,\bar{\beta}(u^\epsilon)})|\ind_{[r\geq t']}\}dr\\ \\ -\integ{0}{\rho_n^\epsilon\wedge T}|f(r,X_{r}^{t,x,u^{\epsilon},\beta^{\epsilon}({u^\epsilon})})|\ind_{[t\leq r<t']}dr+
 \sum\limits_{1\leq l\leq n}\chi(t\vee\rho_l^\epsilon,\eta^\epsilon_{l})\ind_{[\rho_l^\epsilon\leq T]}\\ \\-\sum\limits_{1\leq l\leq n}\chi(t'\vee\rho_l^\epsilon,\eta^\epsilon_{l})\ind_{[\rho_l^\epsilon\leq T]}-\ind_{[\rho_n^\epsilon\leq T]}\{|(V^{-})_*(\rho_n^\epsilon,X_{\rho_n^\epsilon}^{t,x,u^\epsilon,\beta^{\epsilon}(u^\epsilon)})|\\ \\+|(V^{-})^*(\rho_n^\epsilon,X_{\rho_n^\epsilon}^{t',x,u^\epsilon,\bar{\beta}(u^\epsilon)})|\}- \ind_{[\rho_n^\epsilon=+\infty]}|g(X_{T}^{t,x,u^\epsilon,\beta^\epsilon(u^\epsilon)})-g(X_{T}^{t',x,u^\epsilon,\bar{\beta}(u^\epsilon)})| \bigg]-2\epsilon\qquad\qquad
 \\ \\ \geq \mathbb{E}\bigg[
-\integ{0}{\rho_n^\epsilon\wedge T}\displaystyle\{|f(r,X_{r}^{t,x,u^{\epsilon},\beta^{\epsilon}({u^{\epsilon}})})-f(r,X_{r}^{t',x,\bar{u},u^\epsilon,\bar{\beta}(u^\epsilon)})|\ind_{[r\geq t']}\}dr\\ \\-\integ{0}{\rho_n^\epsilon\wedge T}|f(r,X_{r}^{t,x,u^{\epsilon},\beta^{\epsilon}({u^\epsilon})})|\ind_{[t\leq r<t']}dr-
 n\max\limits_{1\leq l\leq n}\{\sup\limits_{s\leq T} |\chi(t\vee s,\eta^\epsilon_{l})-\chi(t'\vee s,\eta^\epsilon_{l})|\}\\ \\ -\ind_{[\rho_n^\epsilon\leq T]}\{|(V^{-})_*(\rho_n^\epsilon,X_{\rho_n^\epsilon}^{t,x,u^\epsilon,\beta^{\epsilon}(u^\epsilon)})|+|(V^{-})^*(\rho_n^\epsilon,X_{\rho_n^\epsilon}^{t',x,u^\epsilon,\bar{\beta}(u^\epsilon)})|\}\\ \\- \ind_{[\rho_n^\epsilon=+\infty]}|g(X_{T}^{t,x,u^\epsilon,\beta^\epsilon(u^\epsilon)})-g(X_{T}^{t',x,u^\epsilon,\bar{\beta}(u^\epsilon)})| \bigg]-2\epsilon.\qquad\qquad
\end{array}
\end{equation*}
We note that, there exists a constant $C > 0$ such that
$$-\mathbb{E}\big[\ind_{[\rho_n^\epsilon\leq T]}\{|(V^{-})_*(\rho_n^\epsilon,X_{\rho_n^\epsilon}^{t,x,u^\epsilon,\beta^{\epsilon}(u^\epsilon)})|+|(V^{-})^*(\rho_n^\epsilon,X_{\rho_n^\epsilon}^{t',x,u^\epsilon,\bar{\beta}(u^\epsilon)})\}\big]\geq -\frac{C}{n}.$$
Also, taking the limit as $t\rw t'$, and using the uniform continuity of $f$, $g$ and $\chi$ to obtain:
$$
\liminf_{t\rw t'}V^{-}(t,x)\geq V^{-}(t',x)-\frac{C}{n}-2\epsilon.$$
As n and $\epsilon$ are arbitrary then sending $\epsilon\rightarrow 0 $ and $ n\rightarrow +\infty ,$ to obtain:
$$
\liminf_{t\rw t'}V^{-}(t,x)\geq V^{-}(t',x).$$
 Therefore $V^-$
is
lower semi-continuous.

Now we show that $V^-$ is upper
semi-continuous.\\
Fix an arbitrary $\epsilon>0$, and  we assume that $t < t'$. \\Pick $u^\epsilon\in\widehat{\mathcal{U}}_{t,T}$ such that
\begin{equation*}
\begin{array}{ll}
V^{-}(t,x)-\epsilon\leq \mathbb{E}\bigg[\integ{0}{\tau_n^\epsilon\wedge T}f(r,X_{r}^{t,x,u^\epsilon,\beta(u^\epsilon)})\ind_{[r\geq t]}dr-
\sum\limits_{1\leq m\leq n} c(t\vee\tau^\epsilon_{m} ,\xi^\epsilon_{m})\ind_{[\tau^\epsilon_{m}\leq T]}\prod_{l\geq 1 }\ind_{\{\tau^\epsilon_{m}\ne\rho_{l}\}}\\ +\sum\limits_{l\geq 1} \chi(t\vee\rho_{l},\eta_{l}) \ind_{[\rho_{l}\leq\tau_n^\epsilon\wedge T]}+\ind_{[\tau_n^\epsilon\leq T]}(V^{-})^*(\tau_n^\epsilon,X_{\tau_n^\epsilon}^{t,x,u^\epsilon,\beta(u^\epsilon)})+\ind_{[\tau_n^\epsilon=+\infty]}g(X_{T}^{t,x,u^\epsilon,\beta(u^\epsilon)})\bigg],
\end{array}
\end{equation*}
where $\beta\in\widehat{\mathcal{B}}_{t,T}$ will be chosen later. On the other hand, pick $\beta^\epsilon$ belongs to $\bar{\mathcal{B}}_{t',T}$ such that
\begin{equation*}
\begin{array}{ll}
V^{-}(t',x)+\epsilon\geq \mathbb{E}\bigg[\integ{0}{\tau_n^\epsilon\wedge T}f(r,X_{r}^{t',x,\bar{u},\beta^\epsilon(\bar{u})})\ind_{[r\geq t']}dr
-\sum\limits_{1\leq m\leq n} c(t'\vee\tau_m^\epsilon,\xi^\epsilon_{m}) \ind_{[\tau_m^\epsilon\leq T]}\prod_{l\geq 1 }\ind_{\{\tau_m^\epsilon\ne\rho_{l}^\epsilon\}}\\ +\sum\limits_{l\geq 1}\chi(t'\vee\rho^\epsilon_{l},\eta^\epsilon_{l}) \ind_{[\rho^\epsilon_{l} \leq \tau_n^\epsilon\wedge T]}+\ind_{[\tau_n^\epsilon\leq T]}(V^{-})_*(\tau_n^\epsilon,X_{\tau_n^\epsilon}^{t',x,\bar{u},\beta^\epsilon(\bar{u})})+\ind_{[\tau_n^\epsilon=+\infty]}g(X_{T}^{t',x,\bar{u},\beta^\epsilon(\bar{u})})\bigg],
\end{array}
\end{equation*} where $$\bar{u}=\sum\limits_{1\leq m\leq n} \xi^\epsilon_{m}\ind_{[t'\vee\tau^\epsilon_{m},T]}\in\mathcal{U}_{t',T}.$$
Therefore
\begin{equation}
\begin{array}{ll}
\label{equ22}
V^{-}(t,x)-V^{-}(t',x)\\ \\ \leq \mathbb{E}\bigg[\integ{0}{\tau_n^\epsilon\wedge T}f(r,X_{r}^{t,x,u^\epsilon,\beta^{\epsilon}(u^\epsilon)})\ind_{[r\geq t]}dr-
\sum\limits_{1\leq m\leq n} c(t\vee\tau^\epsilon_{m},\xi^\epsilon_{m})\ind_{[\tau^\epsilon_{m}\leq T]}\prod_{l\geq 1 }\ind_{\{\tau^\epsilon_{m}\ne\rho^\epsilon_{l}\}}\\ \\ +\ind_{[\tau_n^\epsilon\leq T]}(V^{-})^*(\tau_n^\epsilon,X_{\tau_n^\epsilon}^{t,x,u^\epsilon,\beta^{\epsilon}(u^\epsilon)})+\ind_{[\tau_n^\epsilon=+\infty]}g(X_{T}^{t,x,u^\epsilon,\beta^\epsilon(u^\epsilon)})\bigg]\\ \\ -\mathbb{E}\bigg[\integ{0}{\tau^\epsilon_n}f(r,X_{r}^{t',x,\bar{u},\beta^\epsilon(\bar{u})})\ind_{[r\geq t']}dr-
\sum\limits_{1\leq m\leq n}c(t'\vee\tau_m^\epsilon,\xi^\epsilon_{m})\ind_{[\tau^\epsilon_m\leq T]}\prod_{l\geq 1 }\ind_{\{\tau^\epsilon_{m}\ne\rho^\epsilon_{l}\}}\\ \\+\ind_{[\tau_n^\epsilon\leq T]}(V^{-})_*(\tau_n^\epsilon,X_{\tau_n^\epsilon}^{t',x,\bar{u},\beta^{\epsilon}(\bar{u})})+\ind_{[\tau_n^\epsilon=+\infty]}g(X_{T}^{t',x,\bar{u},\beta^\epsilon(\bar{u})})\bigg]+2\epsilon\qquad\qquad\qquad
\end{array}
\end{equation}
\begin{equation}
\begin{array}{ll}
 \leq \mathbb{E}\bigg[
\integ{0}{\tau_n^\epsilon\wedge T}\displaystyle\{|f(r,X_{r}^{t,x,u^{\epsilon},\beta^{\epsilon}({u^{\epsilon}})})-f(r,X_{r}^{t',x,\bar{u},\beta^\epsilon(\bar{u})})|\ind_{[r\geq t']}\}dr\\ \\+\integ{0}{\tau_n^\epsilon\wedge T}|f(r,X_{r}^{t,x,u^{\epsilon},\beta^{\epsilon}({u^\epsilon})})|\ind_{[t\leq r<t']}dr-
 \sum\limits_{1\leq m\leq n}c(t\vee\tau_m^\epsilon,\xi^\epsilon_{m})\ind_{[\tau_m^\epsilon\leq T]}\prod_{l\geq 1 }\ind_{\{\tau^\epsilon_{m}\ne\rho^\epsilon_{l}\}} \\ \\ +\sum\limits_{1\leq m\leq n}c(t'\vee\tau_m^\epsilon,\xi^\epsilon_{m})\ind_{[\tau_m^\epsilon\leq T]}\prod\limits_{l\geq 1 }\ind_{\{\tau^\epsilon_{m}\ne\rho^\epsilon_{l}\}} + \ind_{[\tau_n^\epsilon=+\infty]}|g(X_{T}^{t,x,u^\epsilon,\beta^\epsilon(u^\epsilon)})-g(X_{T}^{t',x,\bar{u},\beta^\epsilon(\bar{u})})| \\ \\+\ind_{[\tau_n^\epsilon\leq T]}\{|(V^{-})^*(\tau_n^\epsilon,X_{\tau_n^\epsilon}^{t,x,u^\epsilon,\beta^{\epsilon}(u^\epsilon)})|+|(V^{-})_*(\tau_n^\epsilon,X_{\tau_n^\epsilon}^{t',x,\bar{u},\beta^\epsilon(\bar{u})})|\}\bigg]+2\epsilon\qquad\qquad
 \\ \\ \leq \mathbb{E}\bigg[
\integ{0}{\tau_n^\epsilon\wedge T}\displaystyle\{|f(r,X_{r}^{t,x,u^{\epsilon},\beta^{\epsilon}({u^{\epsilon}})})-f(r,X_{r}^{t',x,\bar{u},\beta^\epsilon(\bar{u})})|\ind_{[r\geq t']}\}dr\\ \\+\integ{0}{\tau_n^\epsilon\wedge T}|f(r,X_{r}^{t,x,u^{\epsilon},\beta^{\epsilon}({u^\epsilon})})|\ind_{[t\leq r<t']}dr+
 n\max\limits_{1\leq m\leq n}\{\sup\limits_{s\leq T} |c(t'\vee s,\xi^\epsilon_{m})-c(t\vee s,\xi^\epsilon_{m})|\}\\ \\+ \ind_{[\tau_n^\epsilon=+\infty]}|g(X_{T}^{t,x,u^\epsilon,\beta^\epsilon(u^\epsilon)})-g(X_{T}^{t',x,\bar{u},\beta^\epsilon(\bar{u})})|+\ind_{[\tau_n^\epsilon\leq T]}\{|(V^{-})^*(\tau_n^\epsilon,X_{\tau_n^\epsilon}^{t,x,u^\epsilon,\beta^{\epsilon}(u^\epsilon)})|\\ \\ \\ +|(V^{-})_*(\tau_n^\epsilon,X_{\tau_n^\epsilon}^{t',x,\bar{u},\beta^\epsilon(\bar{u})})|\} \bigg]+2\epsilon.\qquad\qquad
\end{array}
\end{equation}
We note that, there exists a constant $C > 0$ such that
$$\mathbb{E}\big[\ind_{[\tau_n^\epsilon\leq T]}\{|(V^{-})^*(\tau_n^\epsilon,X_{\tau_n^\epsilon}^{t,x,u^\epsilon,\beta^{\epsilon}(u^\epsilon)})|+|(V^{-})_*(\tau_n^\epsilon,X_{\tau_n^\epsilon}^{t',x,\bar{u},\beta^\epsilon(\bar{u})})|\}\big]\leq \frac{C}{n}.$$
Also, taking the limit as $t\rw t'$, and using the uniform continuity of $f$, $g$ and $c$ to obtain:
$$
\limsup_{t\rw t'}V^{-}(t,x)\leq V^{-}(t',x)+\frac{C}{n}+2\epsilon.$$
As n and $\epsilon$ are arbitrary then sending $\epsilon\rightarrow 0 $ and $ n\rightarrow +\infty ,$ to obtain:
$$
\limsup_{t\rw t'}V^{-}(t,x)\leq V^{-}(t',x).$$ Therefore $V^-$
is
upper semi-continuous. We then proved that $V^-$ is continuous in $t$.\qquad$\Box$
\begin{cor}
The lower (resp,. upper)
value functions are continuous on $[0,T]\times\mathbb{R}^n,$
therefore we have $V^-=(V^-)_*=(V^-)^*$ (resp,. $V^+=(V^+)_*=(V^+)^*$). As a consequence, $V^-$ (resp,. $V^+$) satisfies the classical dynamic programming principle:   given $0\leq t\leq s\leq T$, $x\in \mathbb{R}^n$, and each family of $\mathbb{Q}_{[t,T]}-$valued $\mathcal{F}_{t,s}$-stopping times $\{\tau^{u,v},(u,v)\in\widehat{\mathcal{U}}_{t,T}\times\widehat{\mathcal{V}}_{t,T}\}$, we have
\begin{equation}
\begin{array}{ll}
V^{-}(t,x)=\inf\limits_{\beta\in \widehat{\mathcal{B}}_{t,T}}\sup\limits_{u\in\widehat{\mathcal{U}}_{t,T}}\mathbb{E}\bigg[\integ{t}{\tau^{u,\beta}}f(r,X_{r}^{t,x,u,\beta(u)})dr-
\sum\limits_{m\geq 1} c(\tau_{m},\xi_{m}) \ind_{[\tau_{m}\leq \tau^{u,\beta}]}\prod_{l\geq 1}\ind_{\{\tau_{m}\ne\rho_{l}\}}\\\qquad\qquad\qquad\qquad +\sum\limits_{l\geq 1} \chi(\rho_{l},\eta_{l}) \ind_{[\rho_{l} \leq \tau^{u,\beta}]}+V^{-}(\tau^{u,\beta},X_{\tau^{u,\beta}}^{t,x,u,\beta(u)})\bigg],
\end{array}
\end{equation}
and
\begin{equation}
\begin{array}{ll}
V^{+}(t,x)=\sup\limits_{\alpha\in\widehat{\mathcal{A}}_{t,T}}\inf\limits_{v\in {\widehat{\mathcal{V}}_{t,T}}}\mathbb{E}\bigg[\integ{t}{\tau^{\alpha,v}}f(r,X_{r}^{t,x,\alpha(v),v})dr-
\sum\limits_{m\geq 1} c(\tau_{m},\xi_{m}) \ind_{[\tau_{m}\leq \tau^{\alpha,v}]}\prod_{l\geq 1}\ind_{\{\tau_{m}\ne\rho_{l}\}}\\\qquad\qquad\qquad\qquad +\sum\limits_{l\geq 1} \chi(\rho_{l},\eta_{l}) \ind_{[\rho_{l}\leq \tau^{\alpha,v}]}+V^{+}(\tau^{\alpha,v},X_{\tau^{\alpha,v}}^{t,x,\alpha(v),v})\bigg].
\end{array}
\end{equation}
\end{cor}
\section{Hamilton-Jacobi-Bellman-Isaacs Equation}
\no

 In this section we prove that
the two value functions are viscosity solutions of the Hamilton-Jacobi-Bellman-Isaacs equation (\ref{eq:HJBI}) associated to the stochastic differential game. We begin with the following Proposition.
\begin{pro}
\label{pro3}
The lower and
upper value functions satisfy the following properties:\\
for all $t\in[0, T)$ and $x \in \mathbb{R}^n$,
\begin{equation}
\label{Hinf}
(i) \; H^\chi_{inf}V(t,x)\geq V(t,x).\qquad\qquad\qquad\qquad\qquad\qquad\qquad\qquad\qquad\qquad
\end{equation}
\begin{equation}
\label{Hsup}
(ii)\;\text{If}\;\, H^\chi_{inf}V(t,x)> V(t,x),\; \text{then}\; H^c_{sup}V(t,x)\leq V(t,x).\qquad\qquad\qquad\qquad
\end{equation}
\end{pro}
$Proof$.
 We make the proof only for $V^{-}$, the other case being analogous.

(i) Let $\eta \in\mathcal{V}$, for $\beta'(u)=\eta\ind_{[t,T]}+\sum\limits_{l\geq 1}\eta'_{l} \ind_{[\rho'_{l},T]}$ we have
\begin{eqnarray*}
V^{-}(t,x)\leq\sup\limits_{u\in\mathcal{U}_{t,T}}J(t, x,u,\beta'(u)).
\end{eqnarray*}
Choose $\beta(u)=\sum\limits_{l\geq 1}\eta'_{l} \ind_{[\rho'_{l},T]}$.
Then
\begin{eqnarray*}
V^{-}(t,x)\leq \sup\limits_{u\in\mathcal{U}_{t,T}}J(t,x+\eta,u,\beta(u))+\chi(t,\eta),
\end{eqnarray*}
from which we deduce that the following inequality holds:
\begin{eqnarray*}
V^{-}(t,x)\leq \inf\limits_{\eta\in \mathcal{V}}\big[V^-(t,x+\eta)+\chi(t,\eta)\big].\qquad\quad
\end{eqnarray*}

(ii) the proof proceeds by a case distinction, for $\tau=t$ in dynamic programming principle, and the suboptimality of multiple impulses at the
same time (Assumptions  (\ref{cout1}) and (\ref{cout2})) we have
\begin{equation*}
\begin{array}{ll}
V^{-}(t,x)=\inf\limits_{\rho\in\{t,+\infty\},\eta\in\mathcal{V} }\sup\limits_{\tau\in\{t,+\infty\},\xi\in\mathcal{U}}\big[-
c(t,\xi)\ind_{[\tau=t]}\ind_{[\rho=+\infty]}\\ \\ \qquad\qquad+\chi(t,\eta) \ind_{[\rho=t]}+V^{-}(t,x+\xi\ind_{[\tau=t]}\ind_{[\rho=+\infty]}+\eta \ind_{[\rho=t]})\big].
\end{array}
\end{equation*}
As a consequence we have
\begin{equation*}
\begin{array}{ll}
V^{-}(t,x)=\inf\limits_{\rho\in\{t,+\infty\}}\bigg[\big[\inf\limits_{\eta\in\mathcal{V}}\{\chi(t,\eta) +V^{-}(t,x+\eta)\}\big]\ind_{[\rho=t]}+\big[\sup\limits_{\tau\in\{t,+\infty\},\xi\in\mathcal{U}}\{-
c(t,\xi)\ind_{[\tau=t]}\\ \qquad\qquad\qquad\qquad\qquad\qquad +V^{-}(t,x+\xi\ind_{[\tau=t]})\}\big]\ind_{[\rho=+\infty]}\bigg].
\end{array}
\end{equation*}
If $H^\chi_{inf}V(t,x)> V(t,x),$
then
\begin{equation*}
\begin{array}{ll}
V^{-}(t,x)=\sup\limits_{\tau\in\{t,+\infty\},\xi\in\mathcal{U}}\{-
c(t,\xi)\ind_{[\tau=t]}+V^{-}(t,x+\xi\ind_{[\tau=t]})\}.
\end{array}
\end{equation*}
Therefore,
\begin{equation*}
\begin{array}{ll}
V^{-}(t,x)\geq\sup\limits_{\xi\in\mathcal{U}}\{-
c(t,\xi)+V^{-}(t,x+\xi)\}.\qquad\qquad \Box
\end{array}
\end{equation*}

\indent Now we prove that the two value functions satisfy in the viscosity  sense.
\begin{theo}
The lower and upper
value functions are viscosity solutions to the Hamilton-Jacobi-Bellman-Isaacs equation (\ref{eq:HJBI}).
\end{theo}
$Proof.$
We give the proof for the lower value function $V^-$, the other case being analogous.
First, we prove
the supersolution property. Suppose $V^--\phi$ achieves a local minimum in $[t_0,t_0+\delta)\times B(x_0, \delta)$ with $V^-(t_0,x_0)=\phi(t_0,x_0)$.
We have by Proposition \ref{pro3},
$$V^-(t_0,x_0)-H^\chi_{inf}V^-(t_0,x_0)\leq 0.$$
If
$$V^-(t_0,x_0)-H^\chi_{inf}V^-(t_0,x_0)=0,$$
then we are done. Now suppose
$$V^-(t_0,x_0)-H^\chi_{inf}V^-(t_0,x_0)<-\epsilon<0,$$
we prove by contradiction that
\begin{equation}
\begin{array}{ll} -\cfrac{\partial \phi}{\partial t}(t_0,x_0)-\mathcal{L}\phi (t_0,x_0)-f(t_0,x_0) \geq 0.
\end{array}
\end{equation}
Suppose otherwise, i.e.,$\, -\cfrac{\partial \phi}{\partial t}(t_0,x_0)-\mathcal{L}\phi (t_0,x_0)-f(t_0,x_0) < 0$. Then without
loss of generality we can assume that $ -\cfrac{\partial \phi}{\partial t}(t,x)-\mathcal{L}\phi(t,x)-f(t,x) < 0$ on $[t_0,t_0+\delta)\times B(x_0, \delta)$.\\
Define the
stopping time $\tau$ by
$$\tau=\inf\{t\in\mathbb{Q}_{[t_0,T]}:X_t\notin B(x_0, \delta)\times[t_0,t_0+\delta) \}\land T.$$
By Ito's formula
\begin{equation}
\begin{array}{llll}
\label{tt}
\mathbb{E}[\phi(\tau,X_\tau^{t_0,x_0,u_0,v_0})]-\phi(t_0,x_0)=\mathbb{E}\bigg[\integ{t_0}{\tau}\bigg(\cfrac{\partial \phi}{\partial t}+\mathcal{L}\phi \bigg)(r,X_{r}^{t_0,x_0,u_0,v_0})dr\bigg].
\end{array}
\end{equation}
Let $\epsilon_1>0$, using the dynamic programming principle, we deduce the existence of a strategy $\beta^\epsilon\in\widehat{\mathcal{B}}_{t_0,T}$ such that
\begin{equation*}
\begin{array}{llllll}
V^{-}(t_0,x_0)\geq \mathbb{E}\bigg[\integ{t_0}{\tau\wedge\rho_1}f(r,X_{r}^{t_0,x_0,u_0,\beta^\epsilon(u_0)})dr+\chi(\rho_1,\eta_1)\ind_{[\rho_1\leq \tau\wedge\rho_1]}\\ \qquad\qquad\qquad+V^{-}(\tau\wedge\rho_1,X_{\tau\wedge\rho_1}^{t_0,x_0,u_0,\beta^\epsilon(u_0)})\bigg]-\epsilon_1
\\ \qquad\qquad \geq \mathbb{E}\bigg[\integ{t_0}{\tau\wedge\rho_1}f(r,X_{r}^{t_0,x_0,u_0,\beta^\epsilon(u_0)})dr+\ind_{[\rho_1\leq \tau]}\big(\chi(\rho_1,\eta_1)+V^{-}(\rho_1,X_{\rho_1}^{t_0,x_0,u_0,\beta^\epsilon(u_0)})\big)\bigg]\\ \qquad\qquad\qquad+\mathbb{E}\bigg[\ind_{[\rho_1 > \tau]}V^{-}(\tau,X_{\tau}^{t_0,x_0,u_0,\beta^\epsilon(u_0)})\bigg]-\epsilon_1
\\ \qquad\qquad \geq \mathbb{E}\bigg[\integ{t_0}{\tau\wedge\rho_1}f(r,X_{r}^{t_0,x_0,u_0,v_0})dr+\ind_{[\rho_1\leq \tau]} H^\chi_{inf}V^{-}(\rho_1,X_{\rho_1}^{t_0,x_0,u_0,v_0})\bigg]\\ \qquad\qquad\qquad+\mathbb{E}\bigg[\ind_{[\rho_1 > \tau]}V^{-}(\tau,X_{\tau}^{t_0,x_0,u_0,v_0})\bigg]-\epsilon_1
\\ \qquad\qquad \geq \mathbb{E}\bigg[\integ{t_0}{\tau\wedge\rho_1}f(r,X_{r}^{t_0,x_0,u_0,v_0})dr+V^{-}(\rho_1\wedge\tau,X_{\rho_1\wedge\tau}^{t_0,x_0,u_0,v_0})\bigg]+\epsilon.\mathbb{P}(\rho_1\leq\tau)-\epsilon_1.
\end{array}
\end{equation*}
Therefore, without loss of generality, we only need to consider $\beta^\epsilon\in\widehat{\mathcal{B}}_{t_0,T}$ such that $\rho_1 > \tau$. Then
$$
\phi(t_0,x_0)=V^{-}(t_0,x_0)\geq \mathbb{E}\bigg[\integ{t_0}{\tau}f(r,X_{r}^{t_0,x_0,u_0,v_0})dr+ V^{-}(\tau,X_{\tau}^{t_0,x_0,u_0,v_0})\bigg]-\epsilon_1
 $$ $$\geq \mathbb{E}\bigg[\integ{t_0}{\tau}f(r,X_{r}^{t_0,x_0,u_0,v_0})dr+ \phi(\tau,X_{\tau}^{t_0,x_0,u_0,v_0})\bigg]-\epsilon_1.
$$
Then from (\ref{tt}) and sending $\epsilon_1\rightarrow 0 $ we get
$$
 0\leq \mathbb{E}\bigg[\integ{t_0}{\tau}-\bigg(\cfrac{\partial \phi}{\partial t}+\mathcal{L}\phi+f\bigg)(r,X_{r}^{t_0,x_0,u_0,v_0})dr\bigg],
$$
which is a contradiction. Therefore, we must have
$
  -\bigg(\cfrac{\partial \phi}{\partial t}+\mathcal{L}\phi+f\bigg)(t_0,x_0)\geq 0
$.
\\
Thanks to Proposition $\ref{pro3}$, we have
$$H^c_{sup}V^-(t_0,x_0)\leq V^-(t_0,x_0).$$
Therefore,
\begin{equation}
\begin{array}{ll}
max\Big\{ min\Big[-\cfrac{\partial \phi}{\partial t}(t_0,x_0)-\mathcal{L}\phi(t_0,x_0)-f(t_0,x_0),V^-(t_0,x_0)-H^c_{sup} V^-(t_0,x_0)\Big],\\ \\ \qquad\qquad\qquad V^-(t_0,x_0)-H^\chi_{inf} V^-(t_0,x_0) \Big\}\geq 0,
\end{array}
\end{equation}
which is the supersolution property. The subsolution property is proved analogously.\qquad$\Box$

Now we give an equivalent of  Hamilton-Jacobi-Bellman-Isaacs equation
(\ref{eq:HJBI}). We consider the new function
$\Gamma$ given by the classical change of variable $\Gamma(t,x) =
\exp(t)V(t, x)$, for any $t\in[0,T]$ and $x\in \mathbb{R}^n$.\\ A second property is given by the:
\begin{pro}
$V$ is a viscosity solution of (\ref{eq:HJBI}) if and only if
$\Gamma$ is a viscosity solution to the  Hamilton-Jacobi-Bellman-Isaacs equation  in $[0,T)\times \mathbb{R}^n$,
\begin{equation}
\begin{array}{ll}\label{I}
max\Big\{ min\Big[\Gamma(t,x)-\cfrac{\partial \Gamma}{\partial t}(t,x)-\mathcal{L}\Gamma(t,x)-\exp(t)f(t,x)\\ \\,\Gamma(t,x)-\tilde{H}^c_{sup} \Gamma(t,x)\Big],\Gamma(t,x)-\tilde{H}^\chi_{inf} \Gamma(t,x) \Big\}= 0,
\end{array}
\end{equation}
where
$$ \tilde{H}^c_{sup}\Gamma(t,x)=\sup\limits_{\xi\in \mathcal{U}}[\Gamma(t,x+\xi)-\exp(t)c(t,\xi)]$$ $$\tilde{H}^\chi_{inf}\Gamma(t,x)=\inf\limits_{\eta\in \mathcal{V}}[\Gamma(t,x+\eta)+\exp(t)\chi(t,\eta)].$$
The terminal
condition for $\Gamma$ is: $\Gamma(T,x)=\exp(T)g(x)$ in
$\mathbb{R}^n.$
\end{pro}

\section{Uniqueness of the solution of Hamilton-Jacobi-\\Bellman-Isaacs equation}
\no

In this section we deal with the issue of uniqueness of the solution of system (\ref{eq:HJBI}) and to do so. We need to impose
additional assumption of cost functions

(\textbf{H6}) There exists a function $h:[0,T]\to(0,\infty)$ such that for all $t\in[0,T]$,
\begin{equation}
\label{cout1+}
\chi(t,\eta_1+\eta_2)\leq\chi(t,\eta_1)+\chi(t,\eta_2)-h(t),
\end{equation}
 and
 \begin{equation}
 \label{cout2+}
c(t,\xi_1+\eta+\xi_2)\leq c(t,\xi_1)-\chi(t,\eta)+c(t,\xi_2)-h(t)
 \end{equation}
for every $\xi_1,\xi_2\in\mathcal{U}$, and $\eta,\eta_1, \eta_2\in\mathcal{V}$.

\indent To prove uniqueness viscosity solution to the HJBI equation (\ref{eq:HJBI}), we begin with the following  technical lemma
 which also appears in [$\cite{[JM]}$, Lemma $4.4$], [$\cite{[ACS]}$, Lemma $6.1$].
\begin{lem}
\label{lemunicite}
Suppose $U,V:[0,T]\times\mathbb{R}^n\to \mathbb{R}$ are uniformly continuous functions viscosity subsolution and a viscosity
supersolution to the HJBI equation (\ref{eq:HJBI}), respectively. Let $(t_0,x_0)\in[0,T]\times\mathbb{R}^n$, be such that
\begin{equation}
\label{<Hinf}
V(t_0,x_0)\geq H^{\chi}_{inf}V(t_0,x_0)
\end{equation}
or
\begin{equation}
\label{>Hsup}
V(t_0,x_0)<H^{\chi}_{inf}V(t_0,x_0),\qquad U(t_0,x_0)\leq H^{c}_{sup}U(t_0,x_0).
\end{equation}
Then for every $\epsilon>0$, there exists $\bar{x}\in\mathbb{R}^n$ and $\delta>0$ such that
$$U(t_0,x_0)-V(t_0,x_0)\leq U(t_0,\bar{x})-V(t_0,\bar{x})+\epsilon$$
and
\begin{equation}
V(t,x)<H^{\chi}_{inf}V(t,x),\qquad U(t,x)> H^{c}_{sup}U(t,x),
\end{equation}
when $(t,x)\in[(t_0-\delta)\vee0,t_0+\delta]\times \bar{B}(\bar{x},\delta),$ with $t_0+\delta < T.$
\end{lem}
$Proof$:
Fix $\epsilon>0,$ We divide the proof into three steps.\\
\textbf{Step 1 } Let (\ref{<Hinf}) hold. Then for $\alpha_1\in(0,1),$ there exists $\eta_0\in \mathcal{V}$ such that
$$ V(t_0,x_0)\geq V(t_0,x_0+\eta_0)+\chi(t_0,\eta_0)-\alpha_1\epsilon.$$
Then, we observe that
\begin{equation*}
\begin{array}{ll} V(t_0,x_0+\eta_0+\eta)+\chi(t_0,\eta)-V(t_0,x_0+\eta_0)\\ \geq V(t_0,x_0+\eta_0+\eta)+\chi(t_0,\eta)+\chi(t_0,\eta_0)-V(t_0,x_0)-\alpha_1\epsilon\\
\geq \chi(t_0,\eta)+\chi(t_0,\eta_0)-\chi(t_0,\eta_0+\eta)-\alpha_1\epsilon.
\end{array}
\end{equation*}
Thus, by  taking $\alpha_1$ sufficiently small and (\ref{cout1+}), we get
$$H^{\chi}_{inf}V(t_0,x_0+\eta_0)-V(t_0,x_0+\eta_0)>0.$$
If $U(t_0,x_0+\eta_0) > H^c
_{sup}U(t_0,x_0+\eta_0)$, we take $(t_0,\bar{x}):=(t_0,x_0+\eta_0).$ Otherwise at\\
$(t_0, x_0+\eta_0)$ condition (\ref{>Hsup}) holds.\\
On the other hand, it’s easy to check that
 \begin{equation*}
 \begin{array}{ll}
 U(t_0,x_0+\eta_0)-V(t_0,x_0+\eta_0)\geq U(t_0,x_0+\eta_0)+\chi(t_0,\eta_0)-V(t_0,x_0)-\alpha_1\epsilon\\ \qquad\qquad\qquad\qquad\qquad\qquad\quad\geq U(t_0,x_0)-V(t_0,x_0)-\alpha_1\epsilon
\end{array}
\end{equation*}

\textbf{Step 2 } Now, suppose that (\ref{>Hsup}) holds at $(t_0,x_0)$, then for $\alpha_2\in(0,1)$,
there exists $\xi_0\in\mathcal{U}$ such that
$$ U(t_0,x_0)\leq U(t_0,x_0+\xi_0)-c(t_0,\xi_0)+\alpha_2\epsilon.$$
Then, we have
 \begin{equation*}
 \begin{array}{ll}
 U(t_0,x_0+\xi_0)-V(t_0,x_0+\xi_0)\geq U(t_0,x_0)+c(t_0,\xi_0)-V(t_0,x_0+\xi_0)-\alpha_2\epsilon\\ \qquad\qquad\qquad\qquad\qquad\qquad\quad\geq U(t_0,x_0)-V(t_0,x_0)-\alpha_2\epsilon.
\end{array}
\end{equation*}
 By  taking $\alpha_2$ sufficiently small and (\ref{cout2+}), we can show that
$$U(t_0,x_0+\xi_0)> H^{c}_{sup}U(t_0,x_0+\xi_0).$$
If $V(t_0,x_0+\xi_0)< H_{inf}^{\chi}V(t_0,x_0+\xi_0)$, we take $(t_0,\bar{x}):=(t_0,x_0+\xi_0)$. Otherwise we can
proceed as in Step 1 and we find $\eta_0\in\mathcal{V}$ such that
$$ V(t_0,x_0+\xi_0)\geq V(t_0,x_0+\xi_0+\eta_0)+\chi(t_0,\eta_0)-\alpha_2\epsilon,$$
$$U(t_0,x_0+\eta_0+\xi_0)-V(t_0,x_0+\eta_0+\xi_0)\geq  U(t_0,x_0)-V(t_0,x_0)-2\alpha_2\epsilon,$$
and
$$
V(t_0,x_0+\xi_0+\eta_0)<H_{inf}^{\chi}V(t_0,x_0+\xi_0+\eta_0).$$
On the other hand, we have
\begin{equation*}
\begin{array}{ll}
U(t_0,x_0+\xi_0+\eta_0+\xi)-c(t_0,\xi)-U(t_0,x_0+\xi_0+\eta_0)\\ \leq U(t_0,x_0+\xi_0+\eta_0+\xi)-c(t_0,\xi)-c(t_0,\xi_0)-U(t_0,x_0+\eta_0)\\ \leq U(t_0,x_0+\xi_0+\eta_0+\xi)-c(t_0,\xi)-c(t_0,\xi_0)+\chi(t_0,\eta_0)-U(t_0,x_0)\\ \leq c(t_0,\xi_0+\eta_0+\xi)-c(t_0,\xi)-c(t_0,\xi_0)+\chi(t_0,\eta_0).
\end{array}
\end{equation*}
Thus, by (\ref{cout2+}), we get
$$U(t_0,x_0+\eta_0+\xi_0) > H^c
_{sup}U(t_0,x_0+\eta_0+\xi_0).$$
Therefore we take $(t_0,\bar{x}):=(t_0,x_0+\xi_0+\eta_0)$.

\textbf{Step 3 }
We can find $\alpha> 0$ such that
\begin{equation*}
\left\{
\begin{array}{ll}
U(t_0,\bar{x})>U(t_0,\bar{x}+\xi)-c(t_0,\xi)+\alpha\qquad \forall \xi\in\mathcal{U},\\
V(t_0,\bar{x})<V(t_0,\bar{x}+\eta)+\chi(t_0,\eta)-\alpha\qquad \forall \eta\in\mathcal{V}.
\end{array}
\right.
\end{equation*}
Thus, by uniform continuity  of $U, V, c,$ and $\chi$ we have
\begin{equation*} \left\{
\begin{array}{ll}
U(t,x)>U(t,x+\xi)-c(t,\xi)-2u(|x-\bar{x}|\vee |t-t_0|)-h(|t-t_0|)+\alpha,\\
V(t,x)<V(t,x+\eta)+\chi(t,\eta)+2v(|x-\bar{x}|\vee |t-t_0|)+w(|t-t_0|)-\alpha,
\end{array}
\right. \end{equation*}
in which $u,v,h$ and $w$ are the modulus of continuity of $U,V,c$ and $\chi$  respectively.
Therefore, there exists $\delta>0$ such that
\begin{equation}
V(t,x)<H^{\chi}_{inf}V(t,x),\qquad U(t,x)> H^{c}_{sup}U(t,x),
\end{equation}
when $(t,x)\in[(t_0-\delta)\vee0,t_0+\delta]\times \bar{B}(\bar{x},\delta),$ with $t_0+\delta < T.$\qquad$\Box$

We are going now to address the question of uniqueness of the viscosity solution of Hamilton-Jacobi-Bellman-Isaacs equation (\ref{eq:HJBI}). We have the following:
\beth \label{uni1}
The solution in viscosity sense of Hamilton-Jacobi-Bellman-Isaacs equation (\ref{eq:HJBI}) is unique in the space of  bounded uniformly continuous functions  on $[0,T]\times\mathbb{R}^n$.
\eeth {\it Proof}. We will show by contradiction that if $U$ and $V$ is a subsolution and a supersolution respectively for $(\ref{I})$, then $U\leq V.$ Therefore if we have two solutions of $(\ref{I})$ then they are obviously equal.
Actually for some $R>0$ (large enough) suppose there exists
 $(\hat{t},\hat{x})\in[0,T]\times B(0,R) $ such that $ \sup\limits_{t,x}(U(t,x)-V(t,x))=U(\hat{t},\hat{x})-V(\hat{t},\hat{x})>0$.

\textbf{Step 1. }
 Using Lemma \ref{lemunicite}, we can find  $(\bar{t},\bar{x})\in[0,T)\times\mathbb{R}^n$  and $\delta$ such that
\begin{equation}\label{equuni}
\sup\limits_{I\times\bar{B}_\delta(\bar{x})}(U(t,x)-V(t,x))\geq U(\bar{t},\bar{x})-V(\bar{t},\bar{x}) >0,
\end{equation}
and
\begin{equation}\label{equuni2}
V(t,x)<\tilde{H}^\chi_{inf} V(t,x),\qquad U(t,x)>\tilde{H}^c_{sup} U(t,x),
\end{equation}
for all $(t,x)\in I\times\bar{B}(\bar{x},\delta)$, where $I:=[\bar{t}-\delta,\bar{t}+\delta]\subset[0,T).$\\
Let $(t_0,x_0)\in I\times\bar{B}(\bar{x},\delta)$ such that
\begin{equation*}
\sup\limits_{I\times\bar{B}(\bar{x},\delta)}(U(t,x)-V(t,x))=U(t_0,x_0)-V(t_0,x_0)=\eta>0.
\end{equation*}
For a small $\epsilon,\beta, \theta>0$, let us define:
\begin{equation}
\begin{array}{ll}
\label{phi}
\Phi_{\epsilon}(t,x,y)=U(t,x)-V(t,y)-\displaystyle\frac{1}{2\epsilon}|x-y|^{2}
-\theta(|x-x_0|^{4}+|y-x_0|^{4})\\
\qquad\qquad\qquad\qquad-\beta (t-t_0)^2.
\end{array}
\end{equation}
By the boundedness of $U$ and $V$, that there exists a
$(t_\epsilon,x_\epsilon,y_\epsilon)\in I\times\bar{B}(\bar{x},\delta)\times \bar{B}(\bar{x},\delta) $, attaining the maximum of $\Phi_{\epsilon}$ on $I\times\bar{B}(\bar{x},\delta)\times \bar{B}(\bar{x},\delta)$.\\
On the other hand, from $2 \Phi_{\epsilon}(t_\epsilon,x_\epsilon,y_\epsilon)\geq
\Phi_{\epsilon}(t_\epsilon,x_\epsilon,x_\epsilon)+\Phi_{\epsilon}(t_\epsilon,y_\epsilon,y_\epsilon)$,
we have
\begin{equation}
\displaystyle\frac{1}{2\epsilon}|x_\epsilon -y_\epsilon|^{2} \leq
(U(t_\epsilon,x_\epsilon)-U(t_\epsilon,y_\epsilon))+(V(t_\epsilon,x_\epsilon)-V(t_\epsilon,y_\epsilon)),
\end{equation}
and consequently $\displaystyle\frac{1}{2\epsilon}|x_\epsilon -y_\epsilon|^{2}$ is bounded,
and as $\epsilon\rightarrow 0$, $|x_\epsilon -y_\epsilon|\rightarrow 0$. Since
$U$ and $V$ are uniformly continuous on $I\times\bar{B}(\bar{x},\delta)$, then $\displaystyle\frac{1}{2\epsilon}|x_\epsilon
-y_\epsilon|^{2}\rightarrow 0$ as
$\epsilon\rightarrow 0.$\\
Since
 $\Phi_{\epsilon}(t_\epsilon,x_\epsilon,y_\epsilon)\geq
\Phi_{\epsilon}(t_0,x_0,x_0)$, we have
\begin{equation}
U(t_0,x_0)-V(t_0,x_0)\leq
\Phi_{\epsilon}(t_\epsilon,x_\epsilon,y_\epsilon)\leq
U(t_\epsilon,x_\epsilon)-V(t_\epsilon,y_\epsilon),
\end{equation}
it follow from the continuity of $U$ and $V$ that, up to a
subsequence,
\begin{equation}
\begin{array}{ccc}
 \label{subsequence}
 (t_\epsilon,x_\epsilon,y_\epsilon)\rightarrow (t_0,x_0,x_0)\\
 \theta(|x_\epsilon-x_0|^{4}+|y_\epsilon-x_0|^{4})\rightarrow 0\\
U(t_\epsilon,x_\epsilon)-V(t_\epsilon,y_\epsilon)\rightarrow
U(t_0,x_0)-V(t_0,x_0).
 \end{array}
 \end{equation}
Next, since $x_0\in \bar{B}(\bar{x},\delta)$ then, for $\epsilon$ small enough and at least for a subsequence which we still index by $\epsilon$, we obtain
\begin{equation}\label{equuni2}
V(t_\epsilon,y_\epsilon)<\tilde{H}^\chi_{inf} V(t_\epsilon,y_\epsilon),\qquad U(t_\epsilon,x_\epsilon)>\tilde{H}^c_{sup} U(t_\epsilon,x_\epsilon).
\end{equation}

\textbf{Step 2. } We now show that $t_\epsilon<T.$ Actually if $t_\epsilon
=T$ then,
$$
\Phi_{\epsilon}(t_0,x_0,x_0)\leq
\Phi_{\epsilon}(T,x_\epsilon,y_\epsilon),$$ and,
$$
U(t_0,x_0)-V(t_0,x_0)\leq
\exp(T)g(x_\epsilon) -\exp(T)g(y_\epsilon)- \beta
(T-t_0)^2,
$$
since $U(T,x_\epsilon)=\exp(T)g(x_\epsilon)$, $V(T,y_\epsilon)=\exp(T)g(y_\epsilon)$
and $g$ is uniformly continuous on $\bar{B}(\bar{x},\delta)$. Then as
$\epsilon \rightarrow 0$, we have,
$$
\begin{array}{ll}
\eta &\leq - \beta
(T-t_0)^2,
\end{array}
$$
which yields a contradiction and we have $t_\epsilon \in [0,T)$.\\

\textbf{Step 3 } To complete the proof it remains to show
contradiction. Let us denote
\begin{equation}
\varphi_{\epsilon}(t,x,y)=\displaystyle\frac{1}{2\epsilon}|x-y|^{2}
+\theta(|x-x_0|^{4}+|y-x_0|^{4})+\beta
(t-t_0)^2.
\end{equation}
Then we have: \be \left\{
\begin{array}{lllllll}\label{derive}
D_{t}\varphi_{\epsilon}(t,x,y)=2\beta(t-t_0),\\
D_x\varphi_{\epsilon}(t,x,y)= \displaystyle\frac{1}{\epsilon}(x-y) +4\theta(x-x_0)|x-x_0|^{2}, \\
D_y\varphi_{\epsilon}(t,x,y)= -\displaystyle\frac{1}{\epsilon}(x-y) +
4\theta(y-x_0)|y-x_0|^{2}\\
\\
B(t,x,y)=D^2_{x,y}\varphi_\epsilon(t,x,y)=\displaystyle\frac{1}{\epsilon}\begin{pmatrix}
   I & -I \\
   -I & I
\end{pmatrix}+\begin{pmatrix}
   a(x) & 0 \\
   0 & a(y)
\end{pmatrix}\\
\\
\text{with}\;\;a(x)=4\theta|x-x_0|^2I+8\theta(x-x_0)(x-x_0)^*.
\end{array}
\right. \ee
Then applying the result by Crandall et al. (Theorem 8.3, \cite{[MC]}) to the function
$$U(t,x)-V(t,y)-\varphi_\epsilon(t,x,y)$$
at the point $(t_\epsilon,x_\epsilon,y_\epsilon)$, for any $\epsilon_1>0$, we can find $c,d\in\R$ and $X,Y\in S_n,$ such that:
 \be \left\{
\begin{array}{lllllll}\label{derive1}
\big(c,\displaystyle\frac{1}{\epsilon}(x_\epsilon-y_\epsilon) +4\theta(x_\epsilon-x_0)|x_\epsilon-x_0|^{2},X\big)\in J^{2,+}(U(t_\epsilon,x_\epsilon)), \\
\big(-d,\displaystyle\frac{1}{\epsilon}(x_\epsilon-y_\epsilon) -
4\theta(y_\epsilon-x_0)|y_\epsilon-x_0|^{2},Y\big)\in J^{2,-}(V(t_\epsilon,y_\epsilon)),\\
c+d=D_t\varphi_\epsilon(t_\epsilon,x_\epsilon,y_\epsilon)=2\beta(t_\epsilon-t_0)\;\;\text{and finally}\\
\\
-\big(\displaystyle\frac{1}{\epsilon_1}+\|B(t_\epsilon,x_\epsilon,y_\epsilon)\|\big)I\leq \begin{pmatrix}
   X & 0 \\
   0 & -Y
\end{pmatrix}\leq B(t_\epsilon,x_\epsilon,y_\epsilon)+\epsilon_1B(t_\epsilon,x_\epsilon,y_\epsilon)^2.
\end{array}
\right. \ee
 Then by definition of
viscosity solution, we get:
\begin{equation}\begin{array}{lll}\label{vis_sub1}-c+U(t_\epsilon,x_\epsilon)
-\langle\displaystyle\frac{1}{\epsilon}(x_\epsilon-y_\epsilon) +4\theta
(x_\epsilon-x_0)|x_\epsilon-x_0|^{2},\\\qquad\qquad
b(t_\epsilon,x_\epsilon)\rangle-\displaystyle\frac{1}{2}tr[\sigma^*(t_\epsilon,x_\epsilon)X\sigma(t_\epsilon,x_\epsilon)]-\exp(t_\epsilon)f(t_\epsilon,x_\epsilon)\leq0
\end{array}\end{equation}
 and
\begin{equation}\begin{array}{l}\label{vis_sub11}d+V(t_\epsilon,y_\epsilon)
-\langle\displaystyle\frac{1}{\epsilon}(x_\epsilon-y_\epsilon)-4\theta
(y_\epsilon-x_0)|y_\epsilon-x_0|^{2},\\\qquad\qquad\qquad\qquad
b(t_\epsilon,y_\epsilon)\rangle-\displaystyle\frac{1}{2}tr[\sigma^*(t_\epsilon,y_\epsilon)Y\sigma(t_\epsilon,y_\epsilon)]-\exp(t_\epsilon)f(t_\epsilon,y_\epsilon)\geq0,
\end{array}\end{equation}
which implies that:
\begin{equation}
\begin{array}{llllll}
\label{viscder} &-c-d+U(t_\epsilon,x_\epsilon)-V(t_\epsilon,y_\epsilon)\\& \leq
[\langle\displaystyle\frac{1}{\epsilon}(x_\epsilon-y_\epsilon) ,
b(t_\epsilon,x_\epsilon)-b(t_\epsilon,y_\epsilon)\rangle\\&+\langle 4\theta
(x_\epsilon-x_0)|x_\epsilon-x_0|^{2}, b(t_\epsilon,x_\epsilon)\rangle +
\langle4\theta (y_\epsilon-x_0)|y_\epsilon-x_0|^{2},
b(t_\epsilon,y_\epsilon)\rangle\\
&+\displaystyle\frac{1}{2}tr[\sigma^*(t_\epsilon,x_\epsilon)X\sigma(t_\epsilon,x_\epsilon)-\sigma^*(t_\epsilon,y_\epsilon)Y\sigma(t_\epsilon,y_\epsilon)]\\
&
+\exp(t_\epsilon)f(t_\epsilon,x_\epsilon)-\exp(t_\epsilon)f(t_\epsilon,y_\epsilon)].
\end{array}
\end{equation}
But from $(\ref{derive})$ there exist a constant $C > 0$ such that:
$$(\|a(x_\epsilon)\|\lor \|a(y_\epsilon)\|)\leq C\theta.$$
As
$$ B=B(t_\epsilon,x_\epsilon,y_\epsilon)=\displaystyle\frac{1}{\epsilon}\begin{pmatrix}
   I & -I \\
   -I & I
\end{pmatrix}+\begin{pmatrix}
   a(x_\epsilon) & 0 \\
   0 & a(y_\epsilon)
\end{pmatrix}$$
Then
$$B\leq\frac{1}{\epsilon}\begin{pmatrix}
   I & -I \\
   -I & I
\end{pmatrix}+C\theta I.$$
It follows that:
\begin{equation}
B+\epsilon_1B^2\leq \frac{\epsilon+\epsilon_1}{\epsilon^2}\begin{pmatrix}
   I & -I \\
   -I & I
\end{pmatrix}+C\theta I,
\end{equation}
where $C$ which hereafter may change from line to line. Choosing now $\epsilon_1=\epsilon,$ yields the relation
\begin{equation}
\label{equaB}
B+\epsilon_1B^2\leq \frac{2}{\epsilon}\begin{pmatrix}
   I & -I \\
   -I & I
\end{pmatrix}+C\theta I.
\end{equation}
Now, from (\textbf{H1}), (\ref{derive1}) and (\ref{equaB}) we get:
$$\displaystyle\frac{1}{2}tr[\sigma^*(t_\epsilon,x_\epsilon)X\sigma(t_\epsilon,x_\epsilon)-\sigma^*(t_\epsilon,y_\epsilon)Y\sigma(t_\epsilon,y_\epsilon)]\leq \frac{C}{\epsilon}|x_\epsilon-y_\epsilon|^2+C\theta(1+|x_\epsilon|^2+|y_\epsilon|^2).$$
Next
$$
\langle\frac{1}{\epsilon}(x_\epsilon-y_\epsilon),b(t_\epsilon,x_\epsilon)-b(t_\epsilon,y_\epsilon)\rangle
\leq \frac{C}{\epsilon}|x_\epsilon - y_\epsilon|^{2}.$$ And finally,
\begin{eqnarray*}
\langle 4\theta
(x_\epsilon-x_0)|x_\epsilon-x_0|^{2}, b(t_\epsilon,x_\epsilon)\rangle+
\langle 4\theta (y_\epsilon-x_0)|y_\epsilon-x_0|^{2},
b(t_\epsilon,y_\epsilon)\rangle\\ \qquad\qquad\leq
C\theta(1+|x_\epsilon||x_\epsilon-x_0|^{3}+|y_\epsilon||y_\epsilon-x_0|^{3}).
\end{eqnarray*}
So that by plugging into (\ref{viscder}) we obtain:
\begin{equation}
\begin{array}{llllll}
\label{viscder11}
-2\beta(t_\epsilon-t_0)+U(t_\epsilon,x_\epsilon)-V(t_\epsilon,y_\epsilon)\\
\leq \displaystyle\frac{C}{\epsilon}|x_\epsilon-y_\epsilon|^2+C\theta(1+|x_\epsilon|^2+|y_\epsilon|^2)+\displaystyle\frac{C}{\epsilon}|x_\epsilon - y_\epsilon|^{2}\\ \qquad +C\theta(1+|x_\epsilon||x_\epsilon-x_0|^{3}+|y_\epsilon||y_\epsilon-x_0|^{3})
\\ \qquad +\exp(t_\epsilon)f(t_\epsilon,x_\epsilon)-\exp(t_\epsilon)f(t_\epsilon,y_\epsilon).
\end{array}
\end{equation}By sending $\epsilon \rightarrow0$, $\beta \rightarrow0$, $\theta
\rightarrow0$,  and taking into account of the continuity of $f$,
we obtain $\eta \leq 0$, which is a contradiction. The proof of Theorem $\ref{uni1}$ is now complete. \qquad$\Box$
\begin{cor}
The lower and
upper value functions coincide, and the value function of the stochastic differential
game is given by $V(t,x):=V^-(t,x)=V^+(t,x)$ for every $(t,x)\in [0,T)\times\mathbb{R}^n.$
\end{cor}


\begin{thebibliography}{99}
\small
\renewcommand{\baselinestretch}{0.3}
\bibitem{[TB]}T. Bielecki and S. Pliska, Risk sensitive asset management with fixed transaction costs,
Finance Stoch., 4 (2000), pp. 1-33.
\bibitem{[BOI]}B. Bouchard, Introduction to stochastic control of mixed diffusion processes, viscosity solu-
tions and applications in finance and insurance, Lecture Notes Preprint, 2007.
\bibitem{[AC]} A. Cadenillas, T. Choulli, M. Taksar, and L. Zhang, Classical and impulse stochastic
control for the optimization of the dividend and risk policies of an insurance firm, Math.
Finance, 16 (2006), pp. 181-202.
\bibitem{[ACF]} A. Cadenillas and F. Zapatero, Optimal central bank intervention in the foreign exchange
market, J. Econom. Theory, 97 (1999), pp. 218-242.
\bibitem{[YS]}Y-S. A. Chen and X. Guo, impulse control of  multidimensional jump
 diffusions in finite time horison, siam J. control optim.Vol. 51, No. 3, pp. 2638-2663 , (2013).
\bibitem{[ACS]} A. Cosso, Stochastic diferential games involving impulse controls and double-obstacle quasi-
variational inequalities, SIAM J. Control Optim. 51 (2013), no. 3, 2102-2131.
\bibitem{[MC]} M. Crandall, Ishii, H. and P.L. Lions (1992) : Users guide to viscosity solutions of second order partial
differential equations, Bull. Amer. Math. Soc., 27, 1-67.
 \bibitem{[JCI]} J. Cvitanic and I. Karatzas, Backward stochastic differential equations with reflection and
Dynkin games, Ann. Probab., 24 (1996), pp. 2024–2056.
\bibitem{[JE]} J. E. Eastham and K. J. Hastings, Optimal impulse control of portfolios, Math. Oper. Res.,
13 (1988), pp. 588-605.
\bibitem{[BE]} B. El Asri Deterministic minimax impulse control in finite horizon: the viscosity solution approach.
ESAIM: Control Optim. Calc. Var., 19, 63-77.
\bibitem{[NE]} N. El Farouq, G. Barles and P. Bernhard, Deterministic minimax impulse control. Appl. Math. Optim. (2010).
\bibitem{[WH]} W. H. Fleming, P. E. Souganidis, On the existence of value functions of
two-player, zero-sum stochastic differential games, Indiana Univ. Math. J., 38
(1989), no. 2, 293-314.
 \bibitem{[SHM]} S. Hamad\`{e}ne and M. Hassani, BSDEs with two reflecting barriers: The general result, Probab.
Theory Related Fields, 132 (2005), pp. 237–264.
\bibitem{[RI]} R. Isaacs, Differential games. A mathematical theory with applications to warfare
and pursuit, control and optimization, John Wiley and Sons, Inc., New York-
London-Sydney, 1965.
\bibitem{[MJS]} M. Jeanblanc and S. Shiryayev, Optimization of the flow of dividends, Russian Math. Surveys,
50 (1995), pp. 257-277.
\bibitem{[MJP]} M. Jeanblanc-Picqu\'{e}, Impulse control method and exchange rate, Math. Finance, 3 (1993),
pp. 161-177.
\bibitem{[IKJ]} I. Kharroubi, J. Ma, H. Pham, and J. Zhang, Backward SDEs with constrained jumps and
quasi-variational inequalities, Ann. Probab., 38 (2010), pp. 794–840.
\bibitem{[RKP]} R. Korn, Protfolio optimization with strictly positive transaction costs and impulse control,
Finance Stoch., 2 (1998), pp. 85-114.
\bibitem{[RKS]}R. Korn, Some applications of impulse control in mathematical finance, Math.Methods Oper.
Res., 50 (1999), pp. 493-518.
\bibitem{[VL]} V. LyVath, M. Mnif, and H. Pham, A model of optimal portfolio selection under liquidity
risk and price impact, Finance Stoch., (2007), pp. 51-90.
\bibitem{[DC]} D. C. Mauer and A. Triantis, Interactions of corporate financing and investment decisions:
A dynamic framework, J. Finance, 49 (1994), pp. 1253-1277.
\bibitem{[PM]}P. Miguel Almeida Serra Costa Vitoria, A weak dynamic programming principle
for zero-sum stochastic differential games
\bibitem{[AJ]} A. J. Morton and S. Pliska, Optimal portfolio management with fixed transaction costs,
Math. Finance, 5 (1995), pp. 337-356.
\bibitem{[GM]} G. Mundaca and B. Oksendal, Optimal stochastic intervention control with application to
the exchange rate, J. Math. Econom., 29 (1998), pp. 225-243.
\bibitem{[BO]} B. Oksendal and A. Sulem, Optimal consumption and portfolio with both fixed and proportional
transaction costs, SIAM J. Control. Optim., 40 (2002), pp. 1765-1790.
\bibitem{[HP]} H. Pham, Optimal stopping of controlled jump diffusion processes: a viscosity solution approach, Journal of
Mathematical Systems, Estimation, and Control, 8 (1998), p. 27 pp. (electronic).
\bibitem{[RC]}R.C. Seydel, Existence and uniqueness of viscosity solutions for QVI associated
with impulse control of jump-diffusions. Stochastic Process. Appl., 119(10):3719-3748,
2009.
\bibitem{[ST]} S. Tang, J. M. Yong, Finite horizon stochastic optimal switching and impulse
controls with a viscosity solution approach, Stochastics Stochastics Rep.,
45 (1993), no. 3-4, 145-176.
\bibitem{[ATJ]} A. Triantis and J. E. Hodder, Valuing flexibility as a complex option, J. Finance, 45 (1990),
pp. 549-565.
\bibitem{[JM]} J. M. Yong, Zero-sum differential games involving impulse controls, Appl. Math. Optim., 29 (1994), pp. 243-261.
\end{thebibliography}
\end{document}